\documentclass[DIV=15, bibliography=totoc]{scrartcl}  
\usepackage[a4paper,margin=1.5cm,bottom=2.5cm]{geometry}
\usepackage{amsmath,amssymb}
\usepackage{algorithmicx}
\usepackage{algpseudocode} 
\usepackage{algorithm}
\usepackage{graphicx, amsthm, bm, lipsum, mathtools}
\usepackage{xfp}
\usepackage[export]{adjustbox}
\usepackage{textcomp}
\usepackage{mathrsfs, bbold}

\usepackage{multicol}
\usepackage{nomencl}
\makenomenclature

\renewcommand{\nompreamble}{\begin{multicols}{2}}
\renewcommand{\nompostamble}{\end{multicols}}
\usepackage{float}
\usepackage{xcolor} 
\definecolor{lightblue}{rgb}{0,0.5,1.0}
\definecolor{linkblue}{rgb}{0,0.1,0.6}
\definecolor{citegreen}{rgb}{0,0.4,0.0}
\definecolor{linkred}{rgb}{0.8,0,0.005}
\definecolor{mailviolet}{rgb}{0.3,0,0.35}
\definecolor{tumblue}{rgb}{0,0.396,0.741}
\definecolor{darkgreen}{rgb}{0,0.4,0} 
\definecolor{darkbrown}{rgb}{0.5, 0.396, 0.09}

\usepackage[hyperindex,colorlinks=true,linkcolor=linkblue,citecolor=citegreen,urlcolor=mailviolet,filecolor=linkred]{hyperref}
\usepackage{doi}
\usepackage{caption, subcaption}
\usepackage[automark]{scrlayer-scrpage} 

\usepackage{siunitx}
\usepackage[square,comma,nonamebreak,compress,numbers]{natbib}

\usepackage{soulpos}
\usepackage{ulem}
\usepackage{authblk} 
\usepackage{float}

\usepackage[inkscapearea=page]{svg}
\usepackage[noabbrev,capitalize]{cleveref}
\usepackage{etoolbox}
\usepackage{appendix}
\AtBeginEnvironment{appendices}{\crefalias{section}{appendix}}
\usepackage{upgreek}
\usepackage{color}
\usepackage{enumitem}
\usepackage{fix-cm} 

\makeatletter
\def\@font@warning#1{}
\makeatother

\usepackage{pgfplots}
\usepackage{pgfplotstable}
\usepackage{filecontents}
\usepackage{tikz}
\usepackage{tikzscale}
\usepackage{tikz-3dplot}

\usetikzlibrary{spy}
\usetikzlibrary{intersections}
\usetikzlibrary{arrows,shapes}
\usetikzlibrary{spy}
\usetikzlibrary{backgrounds}  
\usetikzlibrary{decorations}
\usetikzlibrary{decorations.markings}
\usetikzlibrary{positioning}
\usetikzlibrary{patterns}
\usetikzlibrary{calc} 
\usetikzlibrary{quotes} 
\usetikzlibrary{external} 
\usetikzlibrary{matrix}
\usetikzlibrary{angles}
\usetikzlibrary{3d}

\pgfplotscreateplotcyclelist{customCycleList}
{%
  {blue, mark=o},
  {red, mark=square},
  {darkgreen, mark=diamond},
  {cyan, mark=diamond},
  {darkbrown, mark=triangle*},
  {black, mark=pentagon},
}
\pgfplotsset{compat=1.18} 

\pgfplotsset{every axis/.append style= {
    cycle list name=customCycleList,
}}

\usepackage{abstract}
\usepackage{placeins}

\newcommand{\publicationDate}{\today}

\date{}
\clearpairofpagestyles
\ofoot*{} 
\ifoot*{\footnotesize\textit{Preprint}\hfill\publicationDate} 
\setkomafont{pagefoot}{\normalfont} 

\crefname{paragraph}{paragraph}{paragraphs}
\Crefname{paragraph}{Paragraph}{Paragraphs}

\normalfont\fontsize{11}{13}\selectfont

\title{Lightweight return-mapping surrogates for multiscale plasticity: a practical guide}

\author[1]{Alireza Daneshyar}
\author[1]{Leon Herrmann}
\author[1,2]{ Stefan Kollmannsberger\thanks{Corresponding author}}
\affil[1]{Data Science in Civil Engineering, Bauhaus-Universität Weimar, Germany}
\affil[2]{Chair of Computational Modeling and Simulation, Technical University of Munich, School of Engineering and Design, Germany}

\begin{document}
\vspace{-1.5cm} \normalem \maketitle \vspace{-1.5cm} \hrule 

\section*{Abstract}
This paper presents a practical guide to building lightweight neural-network surrogates for the plastic return-mapping process in concurrent multiscale (FE\textsuperscript{2}) simulations. Rather than proposing a new architecture, we
show how a deliberately simple feed-forward network, structured to mirror the classical return-mapping update, can replace the prohibitively expensive nested fine-scale solves that dominate the cost of conventional FE\textsuperscript{2}
schemes based on FFT homogenization at the meso-scale. We walk through the full workflow: generating training data from incremental homogenization analyses, constructing a compact yet sufficient dataset, embedding material symmetries
directly into the mapping, and deploying the trained network as a user-defined material subroutine (UMAT) in a standard finite-element solver---enabling widespread use. A sensitivity study examines the model's robustness to data density, increment size, and mesh refinement, and we characterize the regimes in which the surrogate holds and where it breaks down. For the macroscopically isotropic, two-dimensional plane-stress setting considered here, the surrogate reproduces the reference response while reducing the per-analysis cost from hours to seconds with speed-ups up to $30,000$ over standard FE\textsuperscript{2}. The approach extends naturally to three dimensions and to weaker symmetry assumptions, given an appropriate sampling strategy and dataset.


\vspace{0.25cm}
\noindent\textit{Keywords:} 
deep learning; neural network surrogate; data-driven plasticity; plastic return-mapping; concurrent multiscale simulation; FFT-based homogenization
\vspace{-0.4cm}

\section{Introduction} \label{sec:Introduction}
With the rise of artificial intelligence, neural network-based constitutive modeling is emerging as a revolutionary paradigm in continuum mechanics. The future of the field can be envisioned as a vast library of highly accurate neural networks, trained from high-fidelity data, from which an analyst can select the appropriate model for a given simulation. This paradigm could largely render empirical relations obsolete, which are inherently constrained by a limited range of validity and are sometimes difficult or even impossible to generalize to complex scenarios. The shift is therefore from traditional closed-form equations toward data-driven nonlinear function approximators that learn the material response directly from nano-, micro-, or meso-scale models or even from experimental data~\cite{Jang2021, Herrmann2024, Herrmann2025, Roemer2025}.

\subsection{FE² multiscale frameworks}
Many physical problems involve complex interactions and behaviors occurring across multiple length scales. This complexity renders standard phenomenological approaches highly challenging and often insufficient~\cite{Lucarini2021}. A valid analysis, therefore, requires resolving the strong coupling between different length scales within a numerical model, which is generally infeasible with direct numerical simulation. Multiscale methods bridge this gap by enabling tractable predictions of bulk material behavior directly from the physics of the underlying microstructure~\cite{Shin2026}. Concurrent multiscale frameworks implement this coupling through a bidirectional transfer of mechanical state between scales. This approach provides the profound advantage that only the constitutive behavior at the finest scale must be explicitly defined~\cite{Hu2024}. The most established technique within these frameworks is the FE$^2$ approach~\cite{Feyel1999, Daneshyar2020, Hartmann2023}. It usually employs the Finite Element Method (FEM) at both the macro and fine scales. Therein, a Representative Volume Element (RVE), which is a small fine-scale volume whose effective behavior represents the heterogeneous material, is assigned to each macroscopic integration point. In concurrent FE$^2$ modeling, a boundary value problem governing the RVE is solved for each integration point of the macroscopic model to determine the local constitutive response~\cite{Christoff2024, Cui2026}. However, this procedure is computationally prohibitive, as solving such problems at every Gauss point creates an immense numerical burden, particularly when the fine-scale material shows a highly nonlinear response~\cite{Eivazi2024}.

\subsection{FFT-based approaches}
A significant leap forward was made by Moulinec and Suquet in~\cite{Moulinec1994, Moulinec1998} through an efficient alternative based on the Fast Fourier Transform (FFT). It offers a compelling computational advantage over conventional FEM by replacing the solution of a global system of equations with local calculations over the pixels/voxels of a simple grid~\cite{Schneider2021}. This leads to a reduced memory footprint and an algorithmic complexity of $\mathcal{O}(n\log{n})$ compared to that of FEM, which typically scales as $\mathcal{O}(n^2)$ due to the assembly and solution of a system of equations~\cite{Lucarini2021}. Beyond raw speed, its grid-based formulation offers significant practical utility by operating directly on pixel/voxel data, allowing grayscale images or CT scans to be used as inputs without needing any boundary-conforming meshes. Furthermore, its structure inherently satisfies periodic boundary conditions, whereas FEM requires explicit treatment of periodicity. In addition, derivative computations in real space transform to basic algebraic operations in Fourier space. Augmented by the advent of high-performance parallel FFT libraries~\cite{Dalcin2019}, the combined benefits in speed, memory, and convenience render the FFT-based approach superior to the FEM-based one~\cite{Schneider2021}. Building on its core advantages, the method has now been successfully applied to a vast range of problems, including recent works on thin plate structures~\cite{Li2025a}, fracture simulation~\cite{Aranda2025}, topology optimization~\cite{Matsui2025}, micropolar elastoplasticity~\cite{Francis2025}, thermal conductivity~\cite{Gehrig2025}, and damage growth in composite materials~\cite{Li2025b}, among many others. Despite efficiency gains, the fundamental bottleneck remains, as high-resolution fine-scale simulations are required at every Gauss point in each iteration. This renders multiscale modeling for nonlinear materials computationally intensive and explains why such methods are still barely used in industrial applications~\cite{Bishara2023}.

\subsection{Neural network-based methods}
Machine learning enables systems to learn statistical patterns directly from data rather than relying on explicitly defined rules. Its data-driven nature provides a powerful framework for constructing models of complex systems, especially when the underlying physical principles are abstract, unknown, or difficult to describe~\cite{Herrmann2025}. In the context of data-driven constitutive modeling, deep neural networks have become the predominant machine learning approach. As universal function approximators~\cite{hornik_multilayer_1989,cybenko_approximation_1989}, neural networks excel at capturing complex nonlinear relationships and have proven highly effective for demanding regression tasks~\cite{Jang2021}. Constitutive modeling is particularly well-suited for such a data-driven methodology. The idea of learning the strain-stress map with neural networks rather than deriving it through hand-crafted empirical rules dates back to~\cite{ghaboussi_knowledgebased_1991} and gained renewed momentum through the data-driven computing paradigm of~\cite{kirchdoerfer_data-driven_2016}.

\subsubsection{Hybrid approaches}

Elastoplasticity is one of the most general constitutive formulations. Its theoretical foundation is structured, built from an elastic law, a yield function, a flow rule, and a hardening model. This structure allows for hybrid approaches, where each component can be independently represented by either a classical phenomenological model or a data-driven formulation~\cite{Fuhg2025}. This flexibility was exploited in~\cite{Furukawa2004}, where the hardening model was replaced by two dedicated neural networks, one for the isotropic part and another for the kinematic part, while the elastic law, yield function, and flow rule were retained in their classical forms. Those neural networks incrementally update the hardening internal variables so that the model can capture complex combined hardening behavior under uniaxial cyclic loading. \cite{Jones2018} implemented an alternative hybrid strategy, using one neural network to learn the elastic stress-strain relationship and another to model the plastic flow rule. In~\cite{Stoffel2019}, the yield function is kept in its phenomenological form, while a single neural network provides both the plastic flow rule and the kinematic hardening updates. Tracking the inelastic deformation of porous materials was the goal of the hybrid model in~\cite{Settgast2019}. The elastic response follows a linear law, while an initial yield surface, defined from finite element simulations of periodic RVEs, determines whether and to what extent the elastic limit is violated. A neural network, trained on data from the same finite element analyses, predicts the updated stress state using the principal strains as input. A similar approach was followed in~\cite{Settgast2020}, but with neural networks representing both the yield function and the updates of the internal variables. \cite{Jang2021} replaced the stress-integration algorithm for $J_2$ associated plasticity under isotropic hardening by a neural network surrogate. Unloading, however, follows a conventional theory-based procedure. Their surrogate is trained on a dataset generated numerically from the classical return-mapping algorithm in principal stress space. In addition to the mentioned research, there are numerous works exclusively aimed at providing data-driven representations using hybrid approaches, including, but not limited to~\cite{Hartmaier2020, Shoghi2022, Schmidt2022, Fazily2023, Lange2025, Eivazi2026}.

\subsubsection{Multiscale surrogate models}

These neural network surrogates are particularly valuable in multiscale analysis, where they replace the costly fine-scale RVE simulation at each Gauss point of the macroscopic model. \cite{ghavamian_accelerating_2019} proposed an FE$^2$ framework in which a recurrent neural network learns the micro-level response for a strain-softening Perzyna viscoplastic material. To capture the history dependence, the network maps strain sequences to stress sequences, trained on a subspace of possible strain paths that includes non-converged Newton--Raphson iterations. Similarly, \cite{liu_learning-based_2022} map the full deformation history to the homogenized stress. But instead of a recurrent architecture, they compress the strain history via Principle Component Analysis and pass the reduced representation through a feed-forward network. \cite{kalina_fetextrmann_2023} developed a physics-constrained surrogate for the macroscopic RVE response within an FE$^2$ scheme for finite-strain hyperelasticity. \cite{li_multiscale_2020} followed a similar multiscale strategy, training a surrogate at the micro level for an inelastic problem involving $J_2$ plasticity with isotropic hardening. A related but architecturally distinct family are deep material networks, which replace the RVE with a network of simple mechanistic building blocks, e.g., laminates, whose parameters are fit offline to the linear elastic response and reused unchanged for nonlinear and inelastic loading~\cite{liu_deep_2019}, ideal as general homogenization surrogates.

\subsubsection{Path-dependent models}

Path-dependency can be modelled by recurrent neural networks through the hidden state of the network. This approach was proposed by~\cite{mozaffar_deep_2019} to map directly from strain to stress sequences, and explored by~\cite{Heider2020}, who combined a feed-forward network for the path-independent elastic response with a recurrent network for the path-dependent elastoplastic behavior. \cite{ghavamian_accelerating_2019} adopted a similar recurrent architecture in their multiscale framework. \cite{bhattacharya_learning_2023} pursue a related recurrent approach, but with an emphasis on interpretability. Their recurrent neural operator is designed so that learns hidden states correspond to meaningful macroscopic internal variables, building on~\cite{bhattacharya_learning_2023} who prove that non--Markovian homogenized constitutive model can be exactly represented by a finite set of ODE-evolving internal variables. Such recurrent models are particularly valuable when the material law itself is unknown. However, when the goal is computational acceleration of known material laws, as proposed in the present work for multiscale schemes, the standard plasticity internal variables are sufficient to capture path dependence without needing recurrent architectures. By contrast, when the material law is unknown, models are typically enhanced by physics-constrained or physics-augmented formulations~\cite{linden_neural_2023,Fuhg2025,Herrmann2025} to ensure physical reliability of the learned predictions.

Relying on internal variables rather than recurrent networks has the added advantage of naturally handling variable step sizes, which would otherwise need to be addressed through neural ordinary differential equations~\cite{chen_neural_2019} or time-resolution-independent neural operators~\cite{hollenweger_temperature-aware_2026}.

\subsection{Outline of the paper}
The primary contribution of the paper at hand is to replace \emph{the entire} return-mapping components, including the yield function, flow rule, and hardening model, with a single neural network surrogate trained on data generated via FFT-based homogenization. This approach demonstrates the essence of a fully data-driven paradigm for computational plasticity by eliminating the need for explicit analytical definitions in phenomenological constitutive modeling. Accordingly, the remainder of the paper is organized as follows. \Cref{sec:Homogenization} outlines the concurrent multiscale framework. \Cref{sec:Network} details the surrogate model, covering the representative volume element, the neural network architecture, the construction of the training dataset, the implementation of the return-mapping algorithm, and the scalability of the presented approach. \Cref{sec:Benchmark} presents the verification of the surrogate against three standard benchmark problems. \Cref{sec:sensitivity} provides a comprehensive sensitivity analysis of the model with respect to mesh size, step size, stress-angle sampling, the proposed paired-increment strategy, and its generalization to out-of-distribution loading. Finally, \Cref{sec:conclusion} presents some concluding remarks.

\section{Concurrent multiscale framework} \label{sec:Homogenization}
\begin{figure}
    \centering
    \input{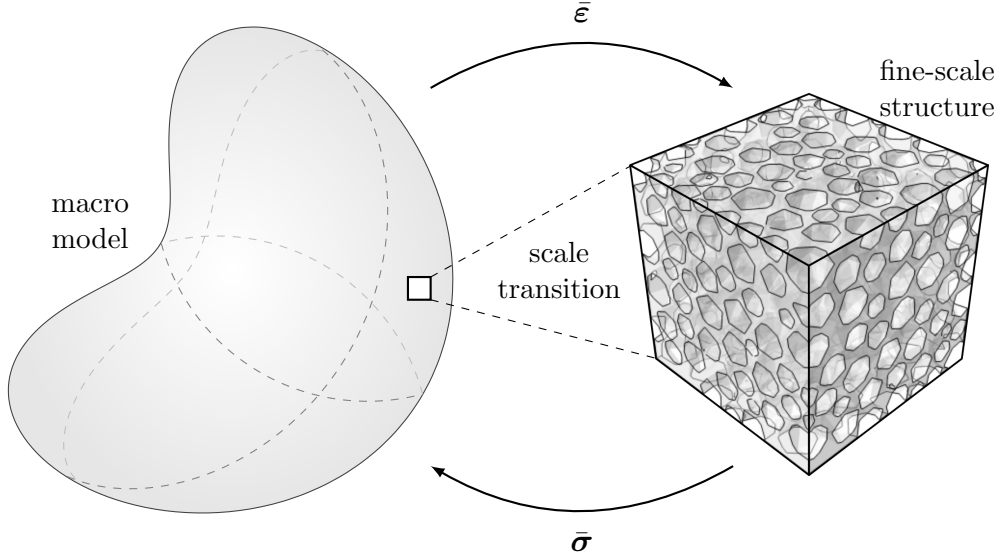}
    \caption{Schematic of the concurrent two-scale coupling: the macro-scale boundary value problem provides the macro strain $\bar{\bm\varepsilon}$ to the fine-scale model, which returns the homogenized stress $\bar{\bm\sigma}$.}
    \label{fig:Concurrent}
\end{figure}

The macroscopic behavior of a material is fundamentally governed by its underlying structure at finer scales. This effect is particularly profound in materials with a relatively coarse microstructure, such as concrete or foam-like materials. Therein, the intrinsic multiscale hierarchy demands a bridge between fine-scale structural details and a continuum-level description. Conventional phenomenological models attempt to establish reliable constitutive relations that approximate the influence of fine-scale features such as grains, fibers, or voids. However, when mechanical phenomena exhibit strong and inseparable coupling across multiple scales, these approximate approaches become inadequate. In such complex scenarios, a concurrent multiscale modeling paradigm becomes indispensable for capturing the integrated material response in an accurate sense. As illustrated in \Cref{fig:Concurrent}, this involves passing the local deformation kinematics from the macro-scale model to a fine-scale simulation and upscaling the resulting homogenized stress---computed by solving the underlying heterogeneous boundary value problem at the fine scale---back to the macro-scale simulation.

\begin{figure}
    \centering
    \input{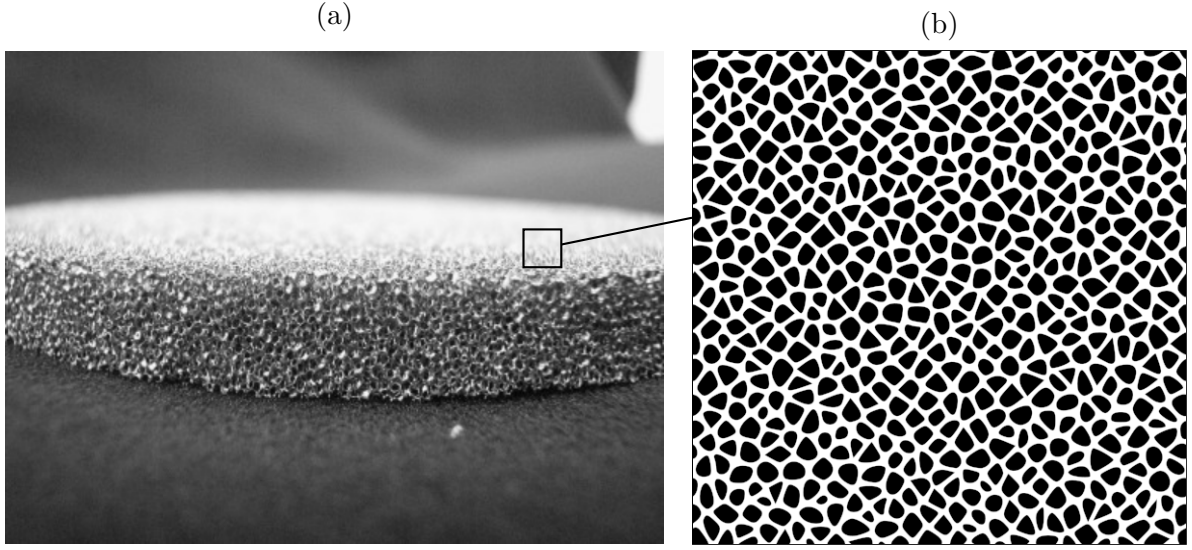}
    \caption{Manufactured cellular metal material: (a) open-cell molybdenum foam~\cite{Stephani2009, Stephani2010} and (b) an example of a two-dimensional RVE with solid phase in white and pore space in black.}
    \label{fig:AluminumFoam}
\end{figure}

Computational homogenization is the standard methodology for this scale transition. It determines the effective macroscopic behavior of a heterogeneous material by solving a boundary value problem of its underlying structure subject to periodic boundary conditions. \Cref{fig:AluminumFoam} shows an RVE of a manufactured open-cell molybdenum foam as an example of a heterogeneous cellular material used in computational homogenization studies. Accurate high-resolution modeling of such fine-scale structures inherently demands significant computational resources, as it requires resolving numerous degrees of freedom and solving complex, often nonlinear material laws. As a result, the practical application of such models depends critically on employing robust and efficient numerical algorithms to deliver reliable solutions within a reasonable simulation time.

To address the above challenge, the initial FEM approach to computational homogenization, such as those used in~\cite{Becker1991, Nakamura1993}, was revolutionized by the FFT-based technique of Moulinec and Suquet~\cite{Moulinec1994, Moulinec1998}. Since then, FFT-based homogenization has become the preferred framework for problems with periodic microstructures~\cite{Schneider2021}, which is addressed in detail in \cref{appendix:fft}.

\section{Neural network surrogate} \label{sec:Network}
We introduce a surrogate modeling strategy to replace concurrent multiscale simulations with a return-mapping-type stress-update procedure. The key departure from conventional approaches is that our method does not require an \textit{a priori} analytical definition of the yield surface or plastic potential. Instead, the underlying constitutive response is learned directly by an artificial neural network trained on data generated via FFT-based homogenization of a chosen RVE.

\subsection{Representative volume element}
To enable a robust and efficient data-driven framework, some primary simplifications are introduced. The first and most critical one is the assumption of material isotropy. For an isotropic material, the mechanical response is invariant with respect to direction; consequently, the yield criterion can be expressed solely in terms of stress invariants. However, constructing an RVE that exhibits a truly isotropic macroscopic response is non-trivial. As a simple counterexample, consider an RVE containing a single, centered circular pore. Although geometrically symmetric in isolation, its periodic arrangement creates a lattice-like structure where the spacing between pores differs along different directions, introducing a directional bias. This anisotropy persists if the RVE contains an oriented pattern of features. Achieving an isotropic response requires a statistically representative RVE such that it must be large enough and contain a sufficiently random, non-periodic distribution of inclusions or pores so that directional dependencies are averaged out. Only such an RVE yields a frame-indifferent, isotropic homogenized response.

\begin{figure}
    \centering
    \begin{tikzpicture}
    \pgfmathsetmacro{\rot}{165}
    \pgfmathsetmacro{\sc}{1.25}
    \pgfmathsetmacro{\dx}{+0.3}
    \pgfmathsetmacro{\dy}{-0.1}

    \begin{scope}[scale=2.0]
        \tdplotsetmaincoords{50}{125}
        \begin{axis}[
            view={125}{50},
            hide axis,
            xmax=2, ymax=2, zmax=2,
            xmin=-2, ymin=-2, zmin=-2,
            colormap={mygray}{color(0)=(black); color(1)=(white)},
            point meta=z,
            point meta min=-3,
            point meta max= 3,
            at={(0,0,0)},
            anchor=origin
        ]
            \addplot3[
                surf,
                samples=2,
                opacity=0.2,
                color=black!20,
                faceted color=black,
                line width=0pt,
                domain=-1:2,
            ] ( {x},{y},{0.25} );
            \addplot3[
                mesh,
                opacity=0.3,
                domain=0:360,
                domain y=0:180,
                line width=0pt,
                mesh/interior colormap name=mygray
            ] ( {(1+0.8*cos(y))*cos(x)*sin(y)}, 
                {(1-0.8*cos(y))*sin(x)*sin(y)}, 
                {1.8*cos(y)} );
            \addplot3[
                smooth,
                color=black,
                domain=0:360,
                domain y=0:180,
                line width=0.25pt,
            ] ( {\sc*cos(\rot)*(0.6*cos(\x)) - \sc*sin(\rot)*((0.8-0.7*sin(\x))*sin(\x)) + \dx},
                {\sc*sin(\rot)*(0.6*cos(\x)) + \sc*cos(\rot)*((0.8-0.7*sin(\x))*sin(\x)) + \dy}, {0.25+0.25*cos(\x)});
        \end{axis}
        \node[font=\normalsize, align=center] at (-1.5,1) {plane stress\\yield limit};
        \draw[draw=black!75, very thin] (-1.25,0.75) -- (-0.75,0);
        \node[font=\normalsize, align=center] at (1.25,1.25) {yield\\locus};

        \draw[->, >=latex, tdplot_main_coords] (0,0,0) -- (2.5,0,0) node[anchor=north east]{$\bar\sigma_x$};
        \draw[->, >=latex, tdplot_main_coords] (0,0,0) -- (0,2.5,0) node[anchor=north west]{$\bar\sigma_y$};
        \draw[->, >=latex, tdplot_main_coords] (0,0,0) -- (0,0,2.0) node[anchor=south]{$\bar\sigma_z$};
    \end{scope}
    \draw[->, >=latex, black] (4,1) to[bend right=-30] (7,1);
    \node[font=\normalsize, align=center] at (5.5,1.75) {$\bar\sigma_z=0$};
    \begin{scope}[scale=2.25, shift={(3.5,-0.5)}]
        \filldraw[fill=black!3]
            plot[
                smooth,
                samples=50,
                domain=0:360
            ] ( {\sc*cos(\rot)*(0.6*cos(\x)) - \sc*sin(\rot)*((0.8-0.7*sin(\x))*sin(\x)) + \dx},
                {\sc*sin(\rot)*(0.6*cos(\x)) + \sc*cos(\rot)*((0.8-0.7*sin(\x))*sin(\x)) + \dy} );
        \draw[->, >=latex] (-0.75,0) -- (1.5,0) node[anchor=north east]{$\bar\sigma_x$};
        \draw[->, >=latex] (0,-0.65) -- (0,2) node[anchor=north east]{$\bar\sigma_y$};
        \node[font=\normalsize, align=center] at (0.7,-0.7) {cross-section};
    \end{scope}
\end{tikzpicture}
    \caption{Visualization of an arbitrary three-dimensional yield surface in Haigh--Westergaard stress space and its planar cross-section at $\bar\sigma_z = 0$ representing its plane stress yield criterion.}
    \label{fig:YieldSurface}
\end{figure}
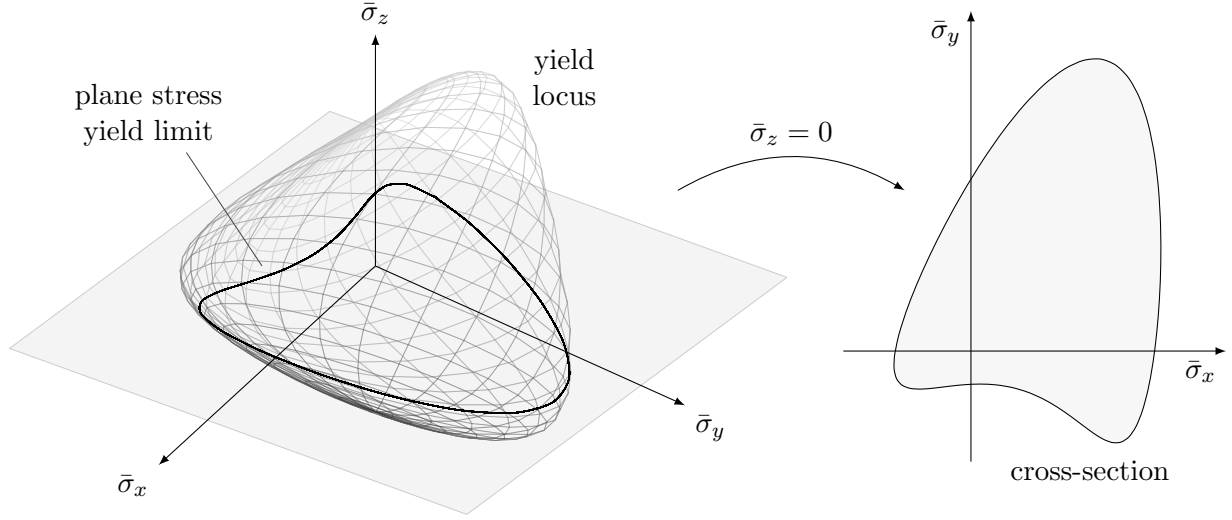

The second simplification is to restrict the analysis to plane stress conditions. This focuses the learning task on a two-dimensional principal stress subspace $(\bar\sigma_x, \bar\sigma_y)$, eliminating the complexity of the out-of-plane stress component and drastically reduces the data required for training. As illustrated in \Cref{fig:YieldSurface}, which shows a three-dimensional yield surface in the Haigh--Westergaard stress space (i.e., the coordinate system having three principal stresses as the axes), the plane-stress yield locus is effectively a cross-section of the full surface at $\bar\sigma_z = 0$. This simplification, however, does not apply to plane strain settings, as the out-of-plane stress remains an active unknown.

An engineering example that obeys these requirements is the RVE of the heterogeneous cellular metal shown in~\Cref{fig:AluminumFoam}b and is investigated in the sequel. It has the topology of a porous matrix (white) embedded in a void phase (black). The solid matrix is modeled as an aluminum alloy with a Young’s modulus $E = 70$ GPa and Poisson’s ratio $\nu = 0.33$, obeying $J_2$ plasticity with a yield stress of 200 MPa and a linear hardening modulus of 700 MPa. Following this characterization, using a $512\times512$ pixel grayscale image, the homogenized stiffness tensor of the RVE expressed in gigapascals reads
\begin{equation}
    \bar{\mathbb{C}} = 
    \left[
    \begin{array}{rrr}
        11.5946 &  7.0232 & 0.0834 \\
         7.0232 & 11.1852 & 0.0840 \\
         0.0834 &  0.0839 & 3.9841
    \end{array}
    \right],
\end{equation}
which deviates from isotropy, yet exhibits near-isotropic behavior. Despite using such an intricate RVE that possesses a well-distributed, sufficiently random arrangement of pores and is large enough to contain a high number of them, the resulting stiffness tensor still shows slight directional dependence. This demonstrates the importance of meticulous RVE design to achieve a truly isotropic macroscopic response. Here, however, we ignore this minor deviation from isotropy and assume a fully isotropic response with no elastic axial-shear coupling in our surrogate. This allows the stress update to be performed in principal stress space without altering the principal directions during the update. For completeness, we note that $\bar{\mathbb{C}}$ is computed through a standard finite differencing procedure: each column corresponds to the solution of an analysis where all strain components are kept zero except for the component corresponding to that column, which is perturbed by an infinitesimally small value. The column entries of $\bar{\mathbb{C}}$ are then obtained by dividing the resulting stress vector by that infinitesimal perturbation.

Beyond the established simplifications, i.e., designing an RVE with sufficiently random features to approximate isotropy and restricting the analysis to the two-dimensional principal stress subspace $(\bar\sigma_x, \bar\sigma_y)$ under plane-stress conditions, a third simplification is the assumption of identical nonlinear response in tension and compression. Although this assumption does not reflect the actual behavior of certain materials such as metallic foams, where compression leads to buckling and collapse of the fine-scale struts, or concrete-like materials and soils that are more sensitive to tension or even cannot withstand it at all, it is reasonable for many others, particularly most bulk metals, which are nearly insensitive to whether the load is tensile or compressive. This assumption of tension--compression symmetry provides a major computational advantage, which will be elaborated in \Cref{subsec:TrainingDataset}.

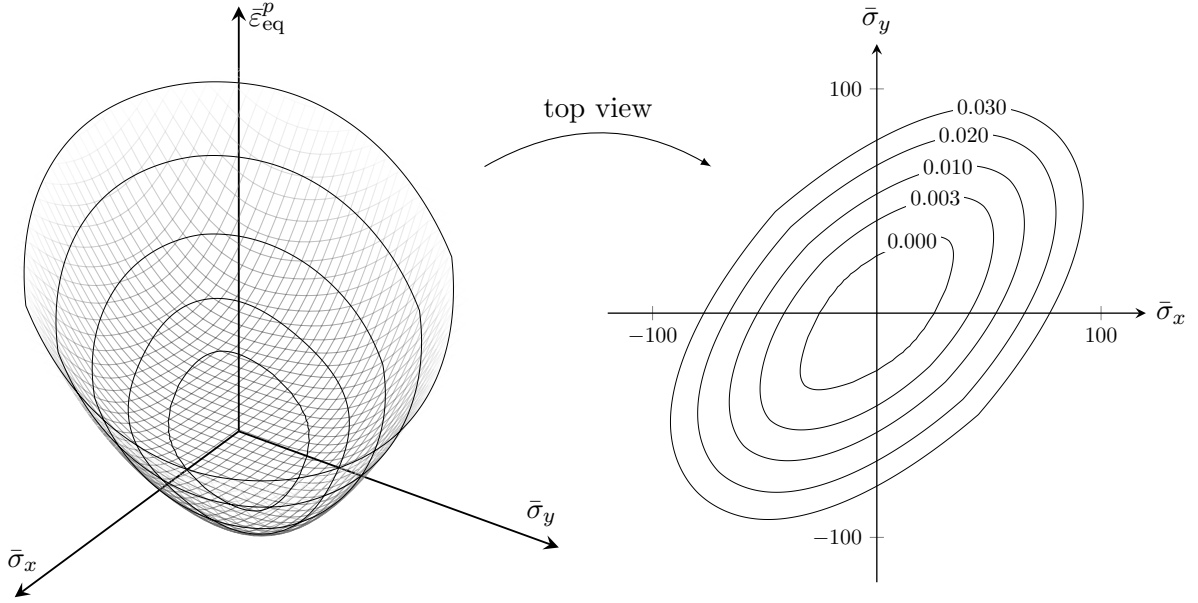
\begin{figure}
    \centering
    \begin{tikzpicture}
    \begin{scope}[scale=1.75]
        \begin{axis}[
            view={125}{50},
            axis lines=middle,
            xmin=0, xmax=110,
            ymin=0, ymax=120,
            zmin=0, zmax=0.085,
            x post scale=0.6,
            y post scale=0.6,
            z post scale=1.5,
            xtick=\empty,
            ytick=\empty,
            ztick=\empty,
            colormap={mygray}{color(0)=(black); color(1)=(white)},
            point meta=z,
            point meta min=0.0,
            point meta max=0.03,
            xlabel={$\bar\sigma_x$},
            ylabel={$\bar\sigma_y$}, 
            zlabel={$\bar\varepsilon^p_{\text{eq}}$},
            xlabel style={yshift=-2pt, scale=0.57},
            ylabel style={yshift=-1pt, scale=0.57},
            zlabel style={yshift= 2pt, scale=0.57},
            anchor=origin,
            unbounded coords=jump,
        ]
            \addplot3[
                mesh,
                opacity=0.25,
                line width=0pt,
                mesh/interior colormap name=mygray
            ]  table {data/evolution-surface.txt};
        \addplot3[smooth, black, line width=0pt] table {data/evolution-levels.txt};
        \end{axis}
    \end{scope}

    \draw[->, >=latex, black] (3.25,3.5) to[bend right=-30] (6.25,3.5);
    \node[align=center] at (4.75,4.25) {top view};

    \begin{scope}[scale=1.25, shift={(6.75,1.25)}]
        \begin{axis}[
            axis lines=middle,
            axis equal image,
            xmin=-120, xmax=120,
            ymin=-120, ymax=120,
            xtick={-100,0,100},
            ytick={-100,0,100},
            ticklabel style={scale=0.6},
            xlabel={$\bar\sigma_x$},
            ylabel={$\bar\sigma_y$},
            xlabel style={anchor=west,  scale=0.8},
            ylabel style={anchor=south, scale=0.8},
            anchor=origin,
            unbounded coords=jump,
        ]
            \addplot[smooth, black, line width=0pt] table {data/evolution-contour.txt};

            \node[fill=white, inner sep=2pt, scale=0.6] at (15.787956, 32.369525) {0.000};
            \node[fill=white, inner sep=2pt, scale=0.6] at (26.158413, 51.338412) {0.003};
            \node[fill=white, inner sep=2pt, scale=0.6] at (31.771317, 65.140396) {0.010};
            \node[fill=white, inner sep=2pt, scale=0.6] at (38.809432, 79.570339) {0.020};
            \node[fill=white, inner sep=2pt, scale=0.6] at (46.877585, 92.002324) {0.030};
        \end{axis}
    \end{scope}    
\end{tikzpicture}
    \caption{Yield surface obtained from RVE analyses and its evolution with accumulated plastic strain.}
    \label{fig:Evolution}
\end{figure}

Based on the chosen RVE, material parameters, and nonlinear constitutive behavior, we construct a three-dimensional representation of the macroscopic yield surface and its evolution due to hardening, as shown in \Cref{fig:Evolution}. To generate this plot, we perform multiple analyses using different combinations of strain inputs applied to the RVE and extract the corresponding principal stresses $\bar\sigma_x$, $\bar\sigma_y$, and the equivalent plastic strain $\bar\varepsilon^p_{\text{eq}}$. The resulting surface is visualized with the principal stresses on the $x$ and $y$ axes and the equivalent plastic strain on the $z$ axis. It should be emphasized that this coordinate system does not correspond to the Haigh--Westergaard stress space, as the $z$ axis represents plastic strain rather than a third principal stress. For additional clarity, a top view of the surface is shown alongside the three-dimensional plot. This view, with axes in megapascals, corresponds to a contour projection of the yield surface onto the stress plane, evolving due to hardening effects. The initial yield surface, corresponding to $\bar\varepsilon^p_{\text{eq}} = 0$, exhibits a pointed oval shape elongated along the line $\bar\sigma_x = \bar\sigma_y$. As plastic strain accumulates, the surface gradually transitions into a smoother, more elliptical form.


There are some aspects that can be concluded from this representation. First of all, although the solid phase of the RVE follows $J_2$ plasticity, the contour lines in the two-dimensional plot deviate from a perfect ellipse, indicating that the macroscopic yield surface renders a material at the global level that is not $J_2$. The second aspect concerns the associativity of the flow rule. Based on the theoretical framework of computational plasticity, the incremental plastic flow rule is given by
\begin{equation} \label{eq:flow}
    \Delta\bar{\bm\varepsilon}^p = \Delta\bar\gamma\frac{\partial\Psi}{\partial\bar{\bm\sigma}},
\end{equation}
where $\Delta\bar\gamma$ is the incremental plastic multiplier and $\Psi$ is a convex function of $\bar{\bm\sigma}$ called the plastic potential~\cite{Neto2011}. The model is considered to be associative if the plastic potential $\Psi$ coincides with the yield function. Otherwise, the plasticity model is considered to be non-associative. Here, the macroscopic yield function and plastic potential are black-box functions, and we do not have analytical representations of them to assess the associativity of the homogenized macroscopic model. Nevertheless, we can analyze this aspect numerically. To this end, we use the homogenized plastic strain increment $\Delta\bar{\bm\varepsilon}^p$ to find the flow direction. However, finding the normal to the yield surface is not as straightforward. The representation in \Cref{fig:Evolution} is merely a geometrical visualization, and those yield contours are extracted from numerous numerical analyses performed on the RVE by collecting stress states that correspond to identical equivalent plastic strain $\bar\varepsilon^p_{\text{eq}}$ to construct those contour lines. Hence, to find the normal to the yield surface, we also adopt a geometrical approach. One way to determine the normal at a given point is to identify two neighboring points on the same contour, one to the left and one to the right, that have exactly the same equivalent plastic strain $\bar\varepsilon^p_{\text{eq}}$. A quadratic curve is then interpolated through these three points, and the normal to this curve at the middle point (the point of interest) provides an approximation of the yield surface normal. 

To compare the evolution of the plastic flow direction with the corresponding normals to the yield surface, we select three representative loading directions. For each case, equal stress increments are applied along a prescribed path in the principal stress plane. To enable interpolation for determining the yield surface normal, additional simulations are performed with perturbations of $\pm 5^\circ$ relative to the main direction. These perturbations provide neighboring points on the same yield contour. For our comparison, we choose three loading directions of $-30^\circ$, $0^\circ$, and $30^\circ$. A fixed stress increment of 20 MPa is used in all cases, and the loading is continued until the equivalent plastic strain $\bar\varepsilon^p_{\text{eq}}$ reaches a value of $0.04$.

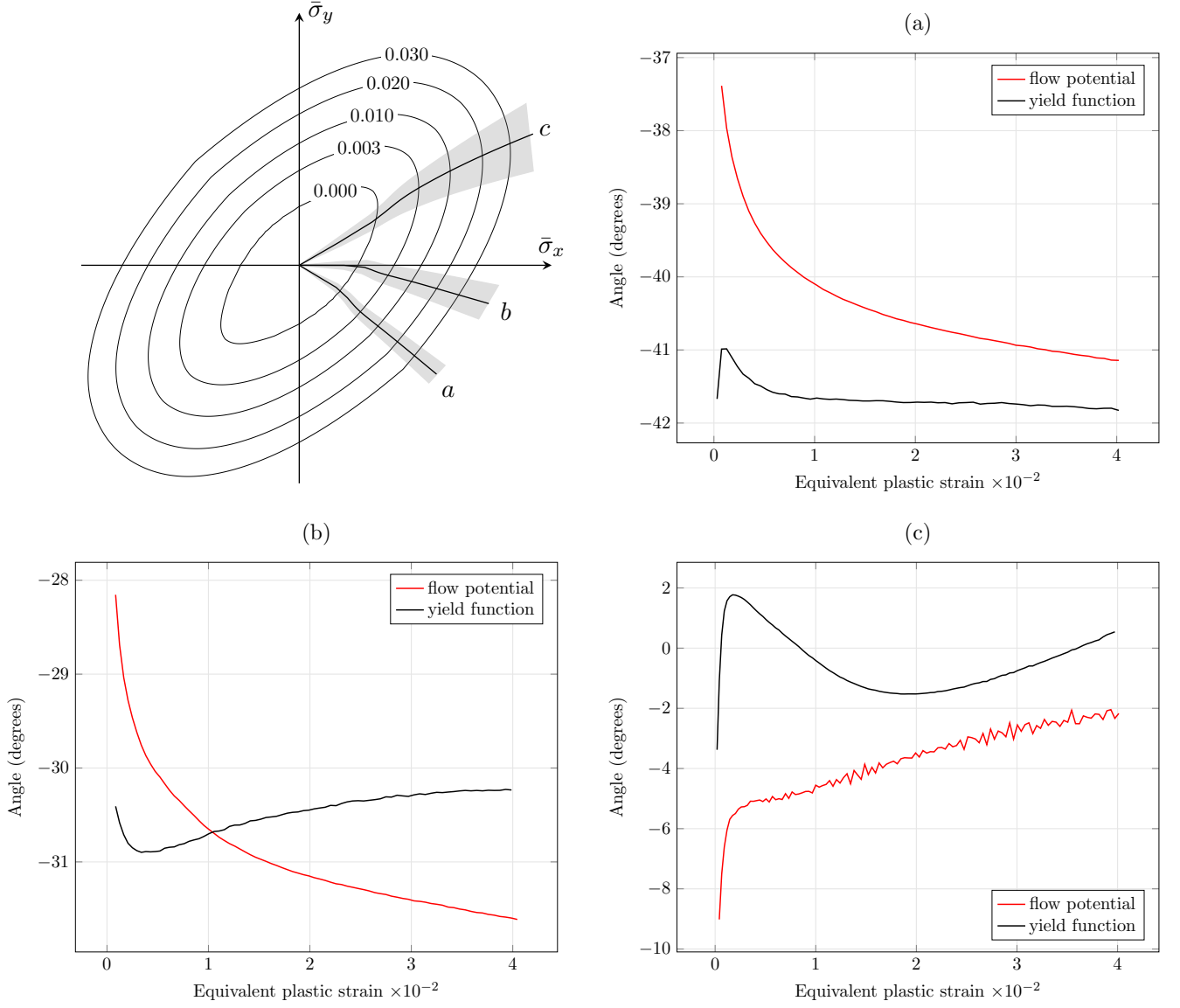
\begin{figure}
    \centering
    \begin{tikzpicture}
    \begin{scope}[scale=1.25, shift={(0,0)}]
        \begin{axis}[
            axis lines=middle,
            axis equal image,
            xmin=-95, xmax=110,
            ymin=-95, ymax=110,
            xtick=\empty,
            ytick=\empty,
            ticklabel style={scale=0.6},
            xlabel={$\bar\sigma_x$},
            ylabel={$\bar\sigma_y$},
            xlabel style={anchor=south,  scale=0.8},
            ylabel style={anchor=west, scale=0.8},
            anchor=center,
            unbounded coords=jump,
        ]
            \addplot[smooth, black, line width=0pt] table {data/evolution-contour.txt};
    
            \node[fill=white, inner sep=2pt, scale=0.6] at (15.787956, 32.369525) {0.000};
            \node[fill=white, inner sep=2pt, scale=0.6] at (26.158413, 51.338412) {0.003};
            \node[fill=white, inner sep=2pt, scale=0.6] at (31.771317, 65.140396) {0.010};
            \node[fill=white, inner sep=2pt, scale=0.6] at (38.809432, 79.570339) {0.020};
            \node[fill=white, inner sep=2pt, scale=0.6] at (46.877585, 92.002324) {0.030};

            \addplot [gray, fill=gray, opacity=0.25, thin, mark=none] table[x index=1, y index=2, col sep=comma] {data/flow-20-n30-adj.txt} -- cycle;
            \addplot [black, thin, mark=none] table[x index=1, y index=2, col sep=comma] {data/flow-20-n30.txt};

            \addplot [gray, fill=gray, opacity=0.25, thin, mark=none] table[x index=1, y index=2, col sep=comma] {data/flow-20-0-adj.txt} -- cycle;
            \addplot [black, thin, mark=none] table[x index=1, y index=2, col sep=comma] {data/flow-20-0.txt};

            \addplot [gray, fill=gray, opacity=0.25, thin, mark=none] table[x index=1, y index=2, col sep=comma] {data/flow-20-30-adj.txt} -- cycle;
            \addplot [black, thin, mark=none] table[x index=1, y index=2, col sep=comma] {data/flow-20-30.txt};

            \node[scale=0.8] (label) at (65, -55) {$a$};
            \node[scale=0.8] (label) at (90, -20) {$b$};
            \node[scale=0.8] (label) at (107, 59) {$c$};

        \end{axis}
    \end{scope}
    
    \begin{scope}[scale=0.7, shift={(13,0)}]
        \begin{axis}[
            width=12cm,
            height=10cm,
            grid=major,
            major grid style={line width=0.6pt, draw=black!10},
            xtick distance=1.0,
            ytick distance=1.0,
            xlabel={Equivalent plastic strain $\times 10^{-2}$},
            ylabel={Angle (degrees)},
            legend pos=north east,
            anchor=center,
            title={(a)},
            title style={scale=1.25},
        ]
            \addplot [red, thick, mark=none] table[x expr=\thisrowno{0}*100, y index=3, col sep=comma, skip first n=3] {data/flow-20-n30.txt};
            \addlegendentry{flow potential}
            
            \addplot [black, thick, mark=none] table[x expr=\thisrowno{0}*100, y index=1, col sep=comma, skip first n=0] {data/flow-20-n30-yield.txt};
            \addlegendentry{yield function}

        \end{axis}
    \end{scope}

    \begin{scope}[scale=0.7, shift={(0,-11)}]
        \begin{axis}[
            width=12cm,
            height=10cm,
            grid=major,
            major grid style={line width=0.6pt, draw=black!10},
            xtick distance=1.0,
            ytick distance=1.0,
            xlabel={Equivalent plastic strain $\times 10^{-2}$},
            ylabel={Angle (degrees)},
            legend pos=north east,
            anchor=center,
            title={(b)},
            title style={scale=1.25},
        ]
            \addplot [red, thick, mark=none] table[x expr=\thisrowno{0}*100, y index=3, col sep=comma, skip first n=4] {data/flow-20-0.txt};
            \addlegendentry{flow potential}
            
            \addplot [black, thick, mark=none] table[x expr=\thisrowno{0}*100, y index=1, col sep=comma, skip first n=2] {data/flow-20-0-yield.txt};
            \addlegendentry{yield function}

        \end{axis}
    \end{scope}    

    \begin{scope}[scale=0.7, shift={(13,-11)}]
        \begin{axis}[
            width=12cm,
            height=10cm,
            grid=major,
            major grid style={line width=0.6pt, draw=black!10},
            xtick distance=1.0,
            ytick distance=2.0,
            xlabel={Equivalent plastic strain $\times 10^{-2}$},
            ylabel={Angle (degrees)},
            legend pos=south east,
            anchor=center,
            title={(c)},
            title style={scale=1.25},
        ]
            \addplot [red, thick, mark=none] table[x expr=\thisrowno{0}*100, y index=3, col sep=comma, skip first n=4] {data/flow-20-30.txt};
            \addlegendentry{flow potential}
            
            \addplot [black, thick, mark=none] table[x expr=\thisrowno{0}*100, y index=1, col sep=comma, skip first n=2] {data/flow-20-30-yield.txt};
            \addlegendentry{yield function}

        \end{axis}
    \end{scope}    
\end{tikzpicture}
    \caption{Associativity analysis of the homogenized model for: (a) $-30^\circ$, (b) $0^\circ$, and (c) $30^\circ$ prescribed loading directions in the principal stress plane.}
    \label{fig:Associativity}
\end{figure}

The resulting stress paths for the three selected loading directions, along with their $\pm5^\circ$ perturbations, are shown in the top left of \Cref{fig:Associativity}. The three main cases correspond to loading directions of $-30^\circ$, $0^\circ$, and $30^\circ$, labeled as (a), (b), and (c), respectively, with their perturbations shaded in gray. Starting from the origin, the stress states initially follow the prescribed linear paths. Once plasticity initiates, they begin to deviate from these directions, resulting in curved trajectories. This occurs under the influence of plastic flow and contrasts sharply with elasticity, where they coincide with the prescribed direction.

Another interesting observation is that the angular spread between the $\pm5^\circ$ perturbations differs across the three cases. For the $-30^\circ$ direction, the spread remains relatively narrow, it becomes wider for $0^\circ$, and even wider for $30^\circ$, with the largest separation between the positive and negative perturbations among all three. This trend indicates that the variation in plastic flow direction is less pronounced around $-30^\circ$, more significant around $0^\circ$, and most pronounced around $30^\circ$.

Finally, we analyze the associativity of the model, i.e., determining whether it exhibits associative or non-associative behavior. To this end, we compare the evolution of the plastic flow direction with the corresponding normals to the yield surface. The results are presented in \Cref{fig:Associativity} for loading directions of $-30^\circ$, $0^\circ$, and $30^\circ$, corresponding to plots (a), (b), and (c), respectively. In each plot, the horizontal axis represents the equivalent plastic strain $\bar\varepsilon^p_{\text{eq}}$, and the vertical axis shows the angle between the corresponding vector and the $\bar\sigma_x$ axis. The red curve tracks the evolution of this angle for the flow direction, while the black curve corresponds to the normal to the yield surface. In all three cases, the flow direction and the yield surface normal do not coincide, nor do they follow similar trends. The trends also differ from one loading path to another. This clear divergence demonstrates that the homogenized macroscopic model exhibits non-associative plastic behavior, despite the underlying model being associative.

As a concluding note, this approach, in which a fixed set of material parameters is used to train a neural network, differs fundamentally from conventional plasticity, where the yield function and plastic potential are analytical functions whose parameters can be adjusted to represent different materials. While the neural network could similarly be trained on a vast dataset encompassing arbitrary combinations of material properties (e.g., varying yield stress and hardening modulus), this generality is arguably unnecessary, especially for combinations that lead to abstract, non-physical material definitions. In practice, engineering materials possess well-defined, fixed properties. A more pragmatic and powerful vision is therefore to develop a comprehensive library of specialized, high-fidelity neural network surrogates, each one representing a specific, real-world material. In this paradigm, an analyst would select a proven, physics-aware model from the library for direct use in their simulation.

\subsection{Neural network architecture}
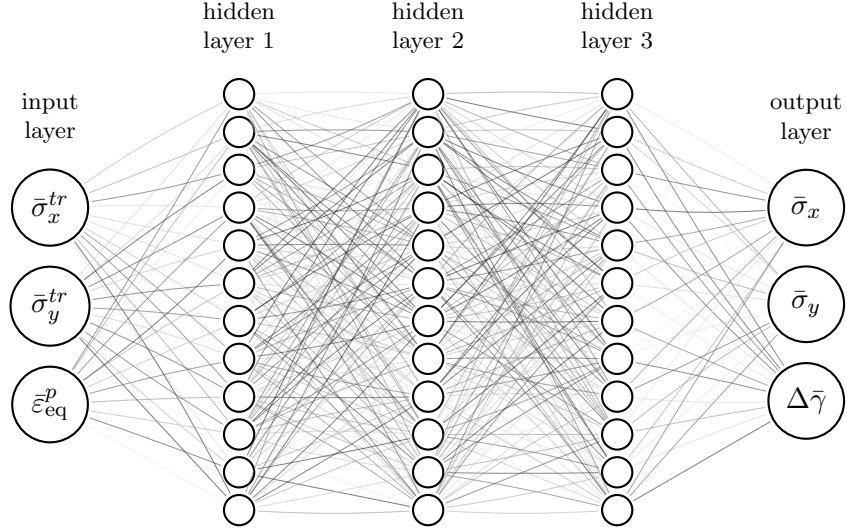
\begin{figure}
    \centering
    \begin{tikzpicture}[
    node distance=1.5cm and 1cm,
    neuron/.style={circle, draw=black, fill=white, minimum size=0.4cm, thick},
    input/.style={circle, draw=black, fill=white, minimum size=1.0cm, thick},
    output/.style={circle, draw=black, fill=white, minimum size=1.0cm, thick},
    layerlabel/.style={text width=2cm, align=center, font=\footnotesize},
]

\pgfmathsetseed{20}

\node[input] (input1) at (0,2.5) {$\bar\sigma^{tr}_x$};
\node[input, below=0.25cm of input1] (input2) {$\bar\sigma^{tr}_y$};
\node[input, below=0.25cm of input2] (input3) {$\bar\varepsilon^p_{\text{eq}}$};
\node[layerlabel, above=0.2cm of input1] (inputlabel) {input\\layer};

\foreach \i in {1,...,12}
    \node[neuron] (h1-\i) at (2.5,{4.5-0.5*\i}) {};
\node[layerlabel, above=0.2cm of h1-1] (hidden1label) {hidden\\layer 1};

\foreach \i in {1,...,12}
    \node[neuron] (h2-\i) at (5,{4.5-0.5*\i}) {};
\node[layerlabel, above=0.2cm of h2-1] (hidden2label) {hidden\\layer 2};

\foreach \i in {1,...,12}
    \node[neuron] (h3-\i) at (7.5,{4.5-0.5*\i}) {};
\node[layerlabel, above=0.2cm of h3-1] (hidden3label) {hidden\\layer 3};

\node[output] (output1) at (10,2.5) {$\bar\sigma_x$};
\node[output, below=0.25cm of output1] (output2) {$\bar\sigma_y$};
\node[output, below=0.25cm of output2] (output3) {$\Delta\bar\gamma$};
\node[layerlabel, above=0.2cm of output1] (outputlabel) {output\\layer};

\foreach \i in {1,2,3}
    \foreach \j in {1,...,12}
        \pgfmathsetmacro{\randcolor}{random(0,100)}
        \pgfmathsetmacro{\randbend}{random(-4,4)}
        \draw[shorten >=1pt, shorten <=1pt,
            bend right=\randbend,
            color=black!\randcolor!white,
            opacity=0.5] (input\i) to (h1-\j);

\foreach \i in {1,...,12}
    \foreach \j in {1,...,12} {
        \pgfmathsetmacro{\randcolor}{random(0,100)}
        \pgfmathsetmacro{\randbend}{random(-4,4)}
        \pgfmathsetmacro{\randdir}{random(0,1)}
        \ifdim\randdir pt<0.5pt
                \draw[shorten >=1pt, shorten <=1pt,
                    bend right=\randbend,
                    color=black!\randcolor!white,
                    opacity=0.5] (h1-\i) to (h2-\j);
        \else \draw[shorten >=1pt, shorten <=1pt,
                    bend left=\randbend,
                    color=black!\randcolor!white,
                    opacity=0.5] (h1-\i) to (h2-\j); \fi}
    
\foreach \i in {1,...,12}
    \foreach \j in {1,...,12} {
        \pgfmathsetmacro{\randcolor}{random(0,100)}
        \pgfmathsetmacro{\randbend}{random(-4,4)}
        \pgfmathsetmacro{\randdir}{random(0,1)}
        \ifdim\randdir pt<0.5pt
                \draw[shorten >=1pt, shorten <=1pt,
                    bend right=\randbend,
                    color=black!\randcolor!white,
                    opacity=0.5] (h2-\i) to (h3-\j);
        \else \draw[shorten >=1pt, shorten <=1pt,
                    bend left=\randbend,
                    color=black!\randcolor!white,
                    opacity=0.5] (h2-\i) to (h3-\j); \fi}

\foreach \i in {1,...,12}
    \foreach \j in {1,2,3}
        \pgfmathsetmacro{\randcolor}{random(0,100)}
        \pgfmathsetmacro{\randbend}{random(-4,4)}
        \draw[shorten >=1pt, shorten <=1pt,
            bend left=\randbend,
            color=black!\randcolor!white,
            opacity=0.5] (h3-\i) to (output\j);

\end{tikzpicture}
    \caption{Neural network architecture with three inputs, three outputs, and three hidden layers having 32 neurons each. For notational clarity, the superscript `\textit{trial}' is shortened to `\textit{tr}' and increment number subscripts are omitted.}
    \label{fig:Architecture}
\end{figure}

The architecture of the neural network is depicted in~\Cref{fig:Architecture} and inspired by the conventional two-step return-mapping procedure in computational plasticity. This standard algorithm employs two possible sets of equations sequentially, selecting the final solution as the only physically valid one. The first step assumes an entirely elastic increment. Consequently, the equivalent plastic strain $\bar\varepsilon^p_{\text{eq}}$ is held fixed, and a trial stress state is computed via
\begin{equation}
    \bar{\bm\sigma}^{trial}_{n+1} = \bar{\bm\sigma}_n + \bar{\mathbb{C}}:\Delta\bar{\bm\varepsilon},
\end{equation}
where subscripts $n$ and $n+1$ denote the increment number and $\Delta\bar{\bm\varepsilon}$ is the macroscopic strain increment. In the conventional scheme, if this elastic trial state satisfies the yield criterion (i.e., if it lies within or on the yield surface), it is accepted as the solution. If it violates the criterion, a plastic corrector step follows. This step provides the updated stress $\bar{\bm\sigma}_{n+1}$ and the incremental plastic multiplier $\Delta\bar\gamma$, so that the plastic consistency condition is satisfied, ensuring the final stress state lies on the updated yield surface.

The presented neural network surrogate adopts, therefore, the identical inputs and outputs as its conventional return-mapping counterpart. The model takes three inputs, two of which are the trial in-plane principal stresses $\bar\sigma^{tr}_x$, $\bar\sigma^{tr}_y$, and the current equivalent plastic strain $\bar\varepsilon^p_{\text{eq}}$. It returns three outputs, comprising the updated in-plane principal stresses $\bar\sigma_x$,$\bar\sigma_y$, and the incremental plastic multiplier $\Delta\bar\gamma$. It is worth noting that in a perfect plasticity scenario, the third input, which is the equivalent plastic strain $\bar\varepsilon^p_{\text{eq}}$, would be redundant as the outputs would not depend on it. Here, however, the presence of hardening renders this internal variable essential. It serves as a state indicator, signaling how far the material has progressed along its hardening path. Conceptually, it corresponds to the expansion of a hypothetical yield surface---a construct that, while not explicitly defined in the presented data-driven model, governs whether a given stress state would cause yielding in the underlying FFT-based fine-scale simulation. The same reasoning applies to the third output, i.e., the incremental plastic multiplier. Under perfect plasticity, its magnitude does not affect the stress solution. Its primary role in that scenario would be to query the spatial distribution and intensity of plastic straining within the RVE.

The presented physics-inspired choice of the input and output layers allows the neural network architecture to provide a very accurate mapping with only three hidden layers of merely 32 neurons. This configuration is sufficiently expressive for highly nonlinear, path-dependent relations, provided that it is trained on large and diverse training datasets. This design choice allows for shifting the focus from architectural optimization of the network itself to the core challenges of this work: investigating the effects of training data generation, performing sensitivity analyses, and establishing robust guidelines for creating comprehensive yet computationally efficient datasets. We use the open-source PyTorch library~\cite{Paszke2019} to train the neural network. To enhance optimization stability, the activation function is chosen to be the Gaussian Error Linear Unit (GELU)~\cite{hendrycks_gaussian_2023} because of its smoothness and differentiability across its entire domain. Furthermore, both inputs and outputs are normalized to compensate for the large dimensional differences between the stress and strain components. This preprocessing prevents the loss landscape from becoming ill-conditioned and ensures that all features are weighted equally during training. The same scaling is retained when integrating the trained surrogate into the finite element solver.

\subsubsection{Physics Constraints}
Embedding physical principles directly into a surrogate reduces the amount of training data required, increases the reliability of its predictions, and has become standard practice in data-driven material modeling~\cite{linden_neural_2023,Fuhg2025,Herrmann2025}. Three main categories of constraints are established:
\begin{enumerate}
    \item \textbf{Invariance}: The response must be unaffected by rigid-body
    transformations of the observer (objectivity, or frame indifference) and by
    the symmetry transformations of the material itself (material symmetry), such
    as isotropy or equal behavior in tension and compression.
    \item \textbf{Thermodynamic consistency}: The response must respect the first
    law (the existence of a free-energy potential $\psi$ from which the stress is
    derived, $\bm{\sigma} = \partial \psi / \partial \bm{\varepsilon}$) and the
    second law (non-negative dissipation, $\mathcal{D} \geq 0$).
    \item \textbf{Material stability}: The response should yield a well-posed
    boundary value problem with a unique solution. This is classically guaranteed
    by convexity~\cite{klein_polyconvex_2022} (of the energy in hyperelasticity, or of the yield surface in
    plasticity), which is desirable for numerical robustness but is not always
    physically accurate.
\end{enumerate}

The simple architecture of \cref{fig:Architecture} already satisfies several of these constraints implicitly. 
Invariance is ensured by mapping from and to ordered principal stresses. Working in the principal frame enforces objectivity: any rotation of the observer's frame leaves the principal stresses $\sigma_1 \geq \sigma_2$ unchanged, since they are computed from the stress tensor's eigenvalues. The principal stress ordering enforces isotropy, since the network cannot tell which spatial direction each eigenvalue originally belonged to, and the sign-based folding of the input domain (to be covered in \cref{subsec:TrainingDataset} and visualized by \cref{fig:StressState}) ensures equal tension--compression behavior. This approach is readily extendable to looser symmetry assumptions, such as anisotropy or distinct tension--compression responses.

Thermodynamic consistency in plasticity requires the plastic dissipation to be non-negative,
\begin{equation}
    \mathcal{D} = \bar{\bm{\sigma}} : \Delta\bar{\bm{\varepsilon}}^{\,p} \geq 0,
    \qquad \Delta\bar{\gamma} \geq 0,
    \label{eq:dissipation}
\end{equation}
of which the irreversibility condition $\Delta\bar{\gamma} \geq 0$ is enforced directly through the fallback to standard elasticity when the network predicts $\Delta\bar{\gamma}< 0$. The dissipation inequality is satisfied only by the training data, but not imposed on the network.\footnote{Importantly, these constraints ensure only physical consistency and not necessarily accuracy for a problem at hand. A simple option, not pursued here, is to use the constraints that are not strongly enforced, i.e., the dissipation inequality and the tension--compression behavior, as indicators of network failure. The network can then report this to the user as a warning.} 

Material stability via convexity of the yield surface is deliberately not enforced. Convexity is overly restrictive for plasticity in homogenized microstructures~\cite{klein_limitations_2026}, whose effective yield surfaces need not be convex, as illustrated by the non-convex homogenized surface in \cref{fig:YieldSurface}.

\subsection{Training dataset} \label{subsec:TrainingDataset}
\begin{figure}
    \centering
    \begin{tikzpicture}
    \begin{scope}[rotate around={45:(0,0)}]
        \fill[gray!30, opacity=0.5, smooth] plot[domain=-90:0, samples=25] 
            ({(2.0 - 0.1*cos(4*\x))*cos(\x)}, {(1.5 + 0.1*cos(4*\x))*sin(\x)}) -- (0,0) -- cycle;
        \draw[thin, black, smooth] plot[domain=0:360, samples=50] 
            ({(2.0 - 0.1*cos(4*\x))*cos(\x)}, {(1.5 + 0.1*cos(4*\x))*sin(\x)});    
            \draw[thin, dashed] (0,-2.5) -- (0,2.5);
            \draw[thin, dashed] (-2.5,0) -- (2.5,0);
            
    \end{scope}
    \draw[->, >=latex] (-3,0) -- (3,0) node[anchor=north east, xshift=1mm]{$\bar\sigma_x$};
    \draw[->, >=latex] (0,-3) -- (0,3) node[anchor=north east, yshift=1mm]{$\bar\sigma_y$};
    \draw[thin] (0,0) -- (2.75,0.5);
    \draw[->] (2.1,0) arc (0:10.3:2.1);
    \node at (2.275,0.2) {$\theta$};
    \filldraw[fill=white, thick] (10.3:1.77) circle (0.075cm);

    \begin{scope}[shift={(3,-3)}, rotate=-45, scale=0.75, local bounding box=mybox]
        \fill[black!5] (-1,-1) rectangle (1,1);
        \draw[thin] (-1,-1) rectangle (1,1);
        \draw[->, very thick, black] (-0.90, 1.25) -- ( 0.90, 1.25); 
        \draw[->, very thick, black] ( 0.90,-1.25) -- (-0.90,-1.25); 
        \draw[->, very thick, black] ( 1.25,-0.90) -- ( 1.25, 0.90); 
        \draw[->, very thick, black] (-1.25, 0.90) -- (-1.25,-0.90); 
    \end{scope}
    \node[font=\tiny, align=center, anchor=north] at (mybox.south) {pure shear};
    
    \begin{scope}[shift={(-3,3)}, rotate=-45, scale=0.75, local bounding box=mybox]
        \fill[black!5] (-1,-1) rectangle (1,1);
        \draw[thin] (-1,-1) rectangle (1,1);
        \draw[->, very thick, black] ( 0.90, 1.25) -- (-0.90, 1.25); 
        \draw[->, very thick, black] (-0.90,-1.25) -- ( 0.90,-1.25); 
        \draw[->, very thick, black] ( 1.25, 0.90) -- ( 1.25,-0.90); 
        \draw[->, very thick, black] (-1.25,-0.90) -- (-1.25, 0.90); 
    \end{scope}
    \node[font=\tiny, align=center, anchor=north] at (mybox.south) {pure shear};
    
    \begin{scope}[shift={(5,0)}, rotate=0, scale=0.75, local bounding box=mybox]
        \fill[black!5] (-1,-1) rectangle (1,1);
        \draw[thin] (-1,-1) rectangle (1,1);
        \draw[->, very thick, black] ( 1.2,0.0) -- ( 2.0,0); 
        \draw[->, very thick, black] (-1.2,0.0) -- (-2.0,0); 
    \end{scope}
    \node[font=\tiny, align=center, anchor=north] at (mybox.south) {uniaxial\\tension};

    \begin{scope}[shift={(-5,0)}, rotate=0, scale=0.75, local bounding box=mybox]
        \fill[black!5] (-1,-1) rectangle (1,1);
        \draw[thin] (-1,-1) rectangle (1,1);
        \draw[->, very thick, black] ( 2.0,0) -- ( 1.2,0.0); 
        \draw[->, very thick, black] (-2.0,0) -- (-1.2,0.0); 
    \end{scope}
    \node[font=\tiny, align=center, anchor=north] at (mybox.south) {uniaxial\\compression};
   
    \begin{scope}[shift={(3,3)}, rotate=45, scale=0.75, local bounding box=mybox]
        \fill[black!5] (-1,-1) rectangle (1,1);
        \draw[thin] (-1,-1) rectangle (1,1);
        \draw[->, very thick, black] (0.0, 1.2) -- ( 0.0, 2.0); 
        \draw[->, very thick, black] (0.0,-1.2) -- (-0.0,-2.0); 
        \draw[->, very thick, black] ( 1.2,0.0) -- ( 2.0, 0.0); 
        \draw[->, very thick, black] (-1.2,0.0) -- (-2.0, 0.0); 
    \end{scope}
    \node[font=\tiny, align=center, anchor=north] at (mybox.south) {equibiaxial\\tension};

    \begin{scope}[shift={(-3,-3)}, rotate=45, scale=0.75, local bounding box=mybox]
        \fill[black!5] (-1,-1) rectangle (1,1);
        \draw[thin] (-1,-1) rectangle (1,1);
        \draw[->, very thick, black] ( 0.0, 2.0) -- (0.0, 1.2); 
        \draw[->, very thick, black] (-0.0,-2.0) -- (0.0,-1.2); 
        \draw[->, very thick, black] ( 2.0, 0.0) -- ( 1.2,0.0); 
        \draw[->, very thick, black] (-2.0, 0.0) -- (-1.2,0.0); 
    \end{scope}    
    \node[font=\tiny, align=center, anchor=north] at (mybox.south) {equibiaxial\\compression};
\end{tikzpicture}
    \caption{Training data domain defined by the stress angle $\theta$, varying from $-45^\circ$ for pure shear stress state to $45^\circ$ for equibiaxial stress state.}
    \label{fig:StressState}
\end{figure}
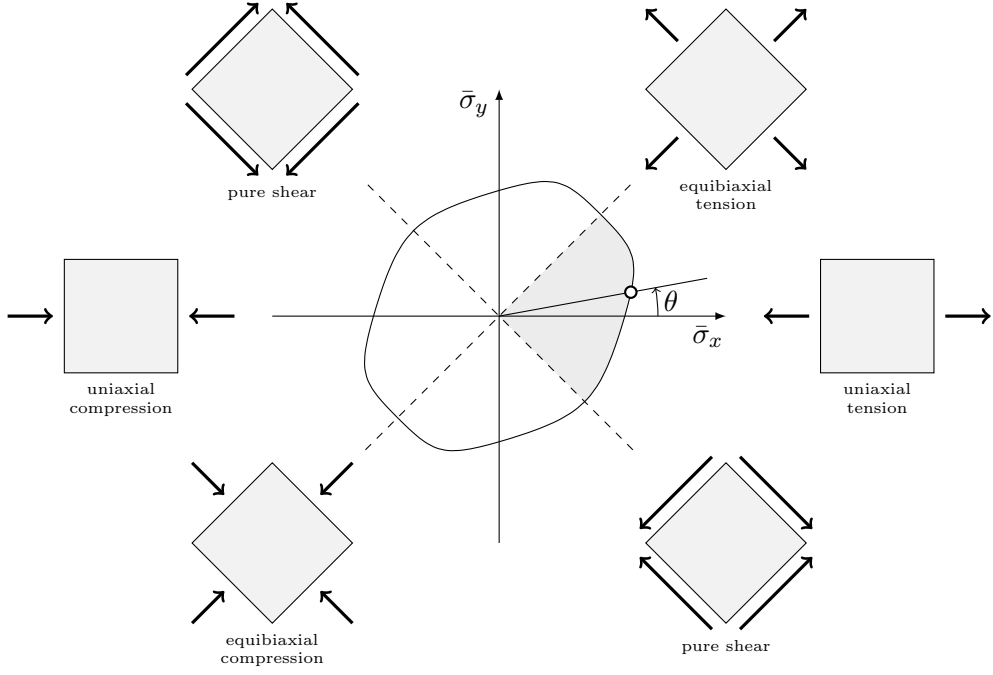

The assumption of identical tensile and compressive behaviors provides a major computational advantage by confining the domain of interest from the full $(\bar\sigma_x, \bar\sigma_y)$ plane to a $90^\circ$ sector. Accordingly, as shown in \Cref{fig:StressState}, the training dataset is generated by varying the stress angle $\theta$ from $-45^\circ$ for pure shear stress state to $45^\circ$ for equibiaxial stress state. Owing to tension--compression symmetry, any arbitrary stress state can be mapped into this canonical sector. This is done by aligning the principal axis corresponding to the maximum absolute principal stress with the $\bar{\sigma}_x$ direction. If this aligned stress is compressive, the stress tensor is multiplied by $-1$ to guarantee that the mapped state resides in the tensile quadrant defined during training. The mapped state is then processed by the trained neural network. Finally, the updated stress is transformed back by re-applying the same sign reversal if necessary. This approach effectively confines the learning domain to a single $90^\circ$ sector while allowing the surrogate model to handle the full stress plane.

We generate the preliminary training dataset by varying the stress angle $\theta$ in $1^\circ$ increments across the defined $90^\circ$ sector. This fine resolution provides a dense data population for the initial verification of the surrogate model against a concurrent multiscale simulation. To investigate the sensitivity of the model's performance to data density, we create additional datasets using coarser angular increments (e.g. $2^\circ$ and $5^\circ$). This sensitivity analysis assesses how the size and distribution of the training data influence the neural network's ability to learn a reliable stress-update mapping.

To generate the training dataset, a series of incremental FFT-based homogenization analyses are performed, each along a predefined stress-angle path defined by a constant $\theta$. This requires imposing the macroscopic strain increment $\Delta\bar{\bm\varepsilon}$ at each loading step. However, the stress update procedure does not inherently preserve a prescribed principal stress ratio. This is a natural aspect of the constitutive response and poses no issue for the dataset itself, where any trial stress state can validly map to its corresponding updated state. Yet, it presents a challenge for data generation along a fixed $\theta$ path. To ensure that each trial data point corresponds to the intended stress direction given the previously deviated updated stress state, the strain increment $\Delta\bar{\bm\varepsilon}$ is adjusted at each step via
\begin{equation}
    \Delta\bar{\bm\varepsilon} = \bar{\mathbb{C}}^{-1}:(\bar{\bm\sigma}^{trial}_{n+1} - \bar{\bm\sigma}_n),
\end{equation}
where the previously updated stress $\bar{\bm\sigma}_n$ generally does not obey the intended stress ratio, yet $\bar{\bm\sigma}^{trial}_{n+1}$ is the target trial stress state aligned with angle $\theta$.

To generate a comprehensive training dataset that enables the surrogate model to accurately handle a wide range of loading increments, a paired-increment strategy is adopted. Each pair consists of an initial, infinitesimally small step that only slightly violates the elastic limit, immediately followed by a substantially larger step that produces noticeable plastic deformation. This approach allows the trained network to accurately interpolate and, to some extent, extrapolate stress updates for arbitrary increment sizes within this range.

It is particularly important to resolve the material response near the onset of plasticity, while larger increments are more efficient once plastic flow is well-developed, avoiding an excessive number of steps to capture advanced plastic deformation. Therefore, the magnitude of the larger step is progressively increased. In this way, the critical elastic-to-plastic transition is sampled densely, while the response under significant yielding is also captured effectively.

A baseline stress increment of 1 MPa is selected empirically. This value remains constant for the small step in each pair throughout the analysis, whereas the large step is computed by augmenting this baseline with an additional part defined by an increment multiplier that grows with the step number. We define this multiplier as
\begin{equation} \label{eq:alpha}
    \alpha = n^2/N,
\end{equation}
where $n$ is the step number and $N$ is the total number of steps. This multiplier progressively increases during the analysis to ensure the desired resolution across the entire loading history. Upon choosing the increment size, it is applied to the larger principal trial stress, $\bar{\sigma}^{tr}_x$, while the smaller principal stress, $\bar{\sigma}^{tr}_y$, is scaled according to the chosen value of $\theta$. An incremental analysis comprising 100 steps is performed along each prescribed $\theta$.

This strategy serves two main purposes. First, it enables the trained network to interpolate stress updates accurately across a range of increment sizes and to generalize robustly beyond the specific loading paths used in training. Second, it provides high resolution near the yield point while efficiently capturing the nonlinear hardening response at large plastic strains, without requiring an excessive number of analyses. Consequently, the resulting dataset spans the full material response, from incipient yielding to fully developed plastic flow.

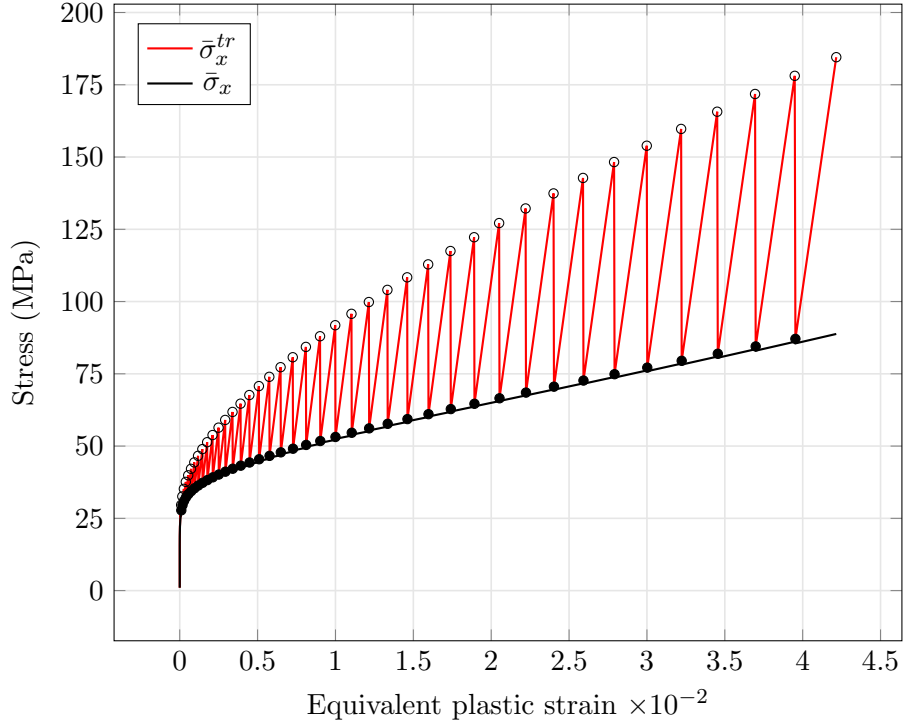
\begin{figure}
    \centering
    \begin{tikzpicture}
    \begin{axis}[
        width=12cm,
        height=10cm,
        grid=major,
        major grid style={line width=0.6pt, draw=black!10},
        xtick distance=0.5,
        ytick distance=25,
        xlabel={Equivalent plastic strain $\times 10^{-2}$},
        ylabel={Stress (MPa)},
        legend pos=north west,
    ]
        \addplot [red, thick, mark=none] table[x expr=(\thisrowno{2}+\thisrowno{5})*100, y index=0, col sep=comma] {data/data-1-0.txt};
        \addlegendentry{$\bar{\sigma}^{tr}_x$}

        \addplot [black, thick, mark=none] table[x expr=(\thisrowno{2}+\thisrowno{5})*100, y index=3, col sep=comma] {data/data-1-0.txt};
        \addlegendentry{$\bar{\sigma}_x$}

        \addplot [black, only marks, mark=o, mark size=1.75pt, mark phase=20, mark repeat=2] table[x expr=(\thisrowno{2}+\thisrowno{5})*100, y index=0, col sep=comma] {data/data-1-0.txt};
        \addplot [black, only marks, mark=*, mark size=1.75pt, mark phase=21, mark repeat=2] table[x expr=(\thisrowno{2}+\thisrowno{5})*100, y index=0, col sep=comma] {data/data-1-0.txt};
    \end{axis}
\end{tikzpicture}
    \caption{Stress versus equivalent plastic strain curves based on the described paired-increment strategy. Filled and hollow circles correspond to the small and large increments, respectively.}
    \label{fig:PairedStrategy}
\end{figure}

\Cref{fig:PairedStrategy} illustrates this paired-increment strategy for $\theta=0^\circ$. The red curve shows the trial stress $\bar{\sigma}^{tr}_x$ as a function of the equivalent plastic strain $\bar\varepsilon^p_{\text{eq}}$, while the black curve shows the updated stress $\bar{\sigma}_x$. The range spanned by the trial and updated stress in each pair enables the model to learn a wide spectrum of strain increments. We later analyze the influence of this small step in the sensitivity analysis section to demonstrate its importance in providing the neural network with the flexibility to handle different strain increments.

Having defined the procedure for adjusting strain increments to align the trial state with the chosen stress angle $\theta$, and having established the details of the paired-increment strategy, the overall process of an incremental FFT-based analysis for generating data along a fixed $\theta$ path is summarized in \Cref{alg:Training}.

\begin{algorithm}
    \begin{algorithmic}[1]
    \State \textbf{Input:} number of steps $N$, baseline increment $\Delta\sigma$, stress angle $\theta$
    \For{$n \in \{0,1,\ldots,N-1\}$}
        \State Define increment size multiplier $\alpha$:
        \If{$n$ is even} 
            \State $\alpha \gets 0$
        \Else
            \State $\alpha \gets n^2/N$
        \EndIf
        \State Compute trial stress along $\theta$:
        \State\qquad $\bar\sigma^{tr}_x \gets \bar\sigma_x + (1+\alpha)\Delta\sigma$
        \State\qquad $\bar\sigma^{tr}_y \gets \bar\sigma^{tr}_x \times\tan\left(\theta\right)$
        \State Compute strain increment $\Delta\bar{\bm{\varepsilon}}$:
        \State\qquad $\Delta\bar{\bm\varepsilon} = \bar{\mathbb{C}}^{-1}:(\bar{\bm\sigma}^{trial}_{n+1} - \bar{\bm\sigma}_n)$
        \State Compute macroscopic strain $\bar{\bm{\varepsilon}}$:
        \State\qquad $\bar{\bm{\varepsilon}} \gets \bar{\bm{\varepsilon}} + \Delta\bar{\bm{\varepsilon}}$
        \State \textbf{goto} \Cref{alg:Procedure} to solve the step
    \EndFor
\end{algorithmic}
    \caption{FFT-based analysis procedure for generating training data along a fixed stress-angle path.}
    \label{alg:Training}
\end{algorithm}

\subsection{Return-mapping algorithm}
The procedure begins by computing the trial stress tensor via the product of the strain increment and the homogenized stiffness tensor, $\bar{\mathbb{C}}$. The principal stresses are then extracted, and the component with the largest absolute value is designated as $\bar{\sigma}_x$, while the other becomes $\bar{\sigma}_y$. To ensure the stress state resides within the canonical sector between $-45^\circ$ and $45^\circ$, the tensor is multiplied by $-1$ if $\bar{\sigma}_x$ is compressive. This mapped stress state, along with the equivalent plastic strain $\bar{\varepsilon}^p_{\text{eq}}$ from the previous converged step, is fed into the trained neural network. The network outputs an updated stress tuple $(\bar{\sigma}_x, \bar{\sigma}_y)$ and the incremental plastic multiplier $\Delta\bar{\gamma}$.

A fundamental postulate of plasticity is the irreversibility of plastic deformation, which mathematically requires the incremental plastic multiplier $\Delta\bar{\gamma}$ to be non-negative. This criterion provides an immediate verification of the prediction provided by the neural network. If the value returned by the network satisfies $\Delta\bar{\gamma} \le 0$, it is interpreted as an indication of purely elastic behavior. In this case, the trial stress state is valid, and the step is considered elastic. Conversely, if $\Delta\bar{\gamma} > 0$, the prediction signals plastic growth, and the corresponding updated stress tuple $(\bar{\sigma}_x, \bar{\sigma}_y)$ is accepted as the correct solution. In this case, the updated stress tuple is transformed back by reversing the initial sign inversion if applied, and then rotated back to the global coordinate system using the principal directions computed prior to the update. The equivalent plastic strain $\bar{\varepsilon}^p_{\text{eq}}$ is also incremented by $\Delta\bar{\gamma}$. Note that unloading would be naturally signaled by a negative $\Delta\bar{\gamma}$. As a result, the surrogate model can track this behavior to some extent. However, capturing complicated unloading and reloading paths requires a designated data generation strategy. This feature is not aimed at in this paper and will be pursued in the future.

It is crucial to emphasize that, as noted earlier, the training dataset enables the network to handle a wide range of strain increments. Nonetheless, accuracy in any return-mapping algorithm, whether conventional or data-driven, requires that loading increments remain within a reasonable bound. To enforce this, we monitor the output $\Delta\bar{\gamma}$. Based on observations from the training data, a practical threshold for this increment can be chosen as $10^{-3}$. If $\Delta\bar{\gamma}$ exceeds this value, the solver is prompted to reject the current load step and re-attempt it with a reduced load increment (halved in our implementation) to ensure stability and preserve solution accuracy.

The return-mapping algorithm presented above is implemented within the Abaqus finite element framework via a user-defined material subroutine (UMAT), written in C++ and released as open-source code. The link to this repository, which includes the implementation and supplementary files, is accessible through \Cref{sec:Data}, the Data Availability section.

\subsection{Scalability} \label{sec:Scalability}
The training dataset employed in this work comprises $91$ stress angles $\times$ $50$ increments $\times$ $2$ (paired) $= 9,100$ data points. To put this number into perspective, a single finite element increment with one global Newton--Raphson iteration requires as many RVE evaluations as there are Gauss points, which amounts to $9,100$ evaluations for a mesh of only $2,275$ elements (${\sim}48\times 48$ grid) with 4 Gauss points per element. Note that this matching of numbers, which corresponds to quite a coarse mesh for a typical finite element analysis, holds only if a one-step analysis is intended and the step is entirely linear over the whole domain (i.e., across all RVEs), such that no more than one global iteration is required. This data efficiency arises directly from the simplifications introduced to demonstrate the viability of simple neural network architectures for accelerating FE$^2$ simulations. We now consider how relaxing each simplification affects the required number of training data points.

Dropping the assumption of identical tension and compression behavior is the most straightforward extension. In that case, we cannot perform the sign reversal to align the maximum absolute principal stress with the $\bar{\sigma}_x$ direction, yet we have the free will to choose which principal stress is along $\bar{\sigma}_x$ and which is along $\bar{\sigma}_y$. Three possibilities arise: one is that both principal stresses are negative. We thus choose the larger one in magnitude to be along $\bar{\sigma}_y$ so that the stress angle $\theta$ falls between $-135^\circ$ and $-90^\circ$. The second case is that one is negative and one is positive, where we choose the positive one to be along $\bar{\sigma}_x$ so that $\theta$ lies in the fourth quadrant. Lastly, if both are positive, we choose the larger one to be along $\bar{\sigma}_x$, hence the stress angle $\theta$ ranges between $0^\circ$ and $45^\circ$. Consequently, we have a symmetry about the $\bar{\sigma}_x = \bar{\sigma}_y$ line, and the stress angle $\theta$ varies between $-135^\circ$ and $45^\circ$ in total. This symmetry is exactly what appears in yield surfaces for pressure-sensitive materials like Drucker--Prager, where the trace on the plane stress surface $\bar{\sigma}_z = 0$ is different in the first and third quadrants, yet is symmetric about the $45^\circ$ line. As a result, dropping the tension-compression symmetry expands the training domain from a $90^\circ$ sector to a half-plane with $180^\circ$ stress angles, increasing the number of data points to $18,100$. This remains well within practical limits and does not alter the fundamental scale of the dataset.

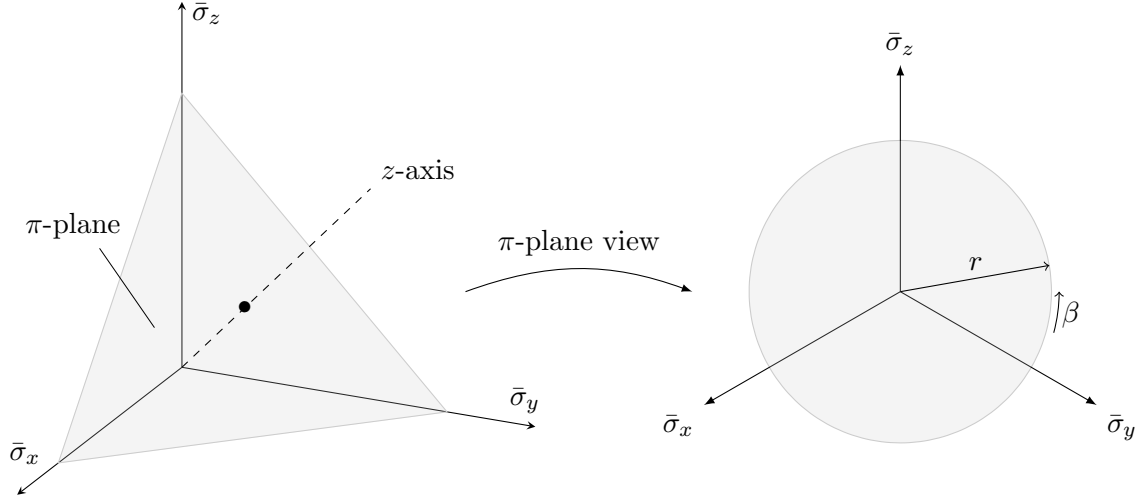
\begin{figure}
    \centering
    \begin{tikzpicture}
    \begin{scope}[scale=1.0]
        \begin{axis}[
            view={115}{30},
            axis lines=middle,
            xmin=0, xmax=3,
            ymin=0, ymax=3,
            zmin=0, zmax=3,
            x post scale=1.0,
            y post scale=1.0,
            z post scale=1.5,
            xtick=\empty,
            ytick=\empty,
            ztick=\empty,
            xlabel={$\bar\sigma_x$},
            ylabel={$\bar\sigma_y$}, 
            zlabel={$\bar\sigma_z$},
            xlabel style={yshift=-2pt, scale=1.0},
            ylabel style={yshift=-1pt, scale=1.0},
            zlabel style={yshift= 2pt, scale=1.0},
            anchor=origin,
        ]
            \addplot3[fill, opacity=0.2, color=black!20, faceted color=black]
                (2.25,0,0) -- (0,2.25,0) -- (0,0,2.25) -- cycle ;
            \addplot3[smooth, color=black!20, very thin]
                (2.25,0,0) -- (0,2.25,0) -- (0,0,2.25) -- cycle ;
            \addplot3[smooth, color=black, dashed]
                (0,0,0) -- (3,3,3);
            \node[scale=1.0] at (1, 1, 1) {$\bullet$};
            \node[scale=1.0, anchor=south west] at (3, 3, 3) {$z$-axis};
            \draw[very thin] (0.5,0,0.5) -- (1.5,0,1.5) node[scale=1.0, anchor=south, xshift=-10]{$\pi$-plane};
        \end{axis}
    \end{scope}
    \draw[->, >=latex, black] (3.75,1) to[bend right=-20] (6.75,1);
    \node[scale=1.0, align=center] at (5.25,1.65) {$\pi$-plane view};
    \begin{scope}[scale=1.0, shift={(9.5, 1.0)}]
        \draw[->, >=latex] (0,0) -- (-3*0.866,-3*0.5) node[anchor=north east]{$\bar\sigma_x$};
        \draw[->, >=latex] (0,0) -- ( 3*0.866,-3*0.5) node[anchor=north west]{$\bar\sigma_y$};
        \draw[->, >=latex] (0,0) -- ( 3*0.0,   3*1.0) node[anchor=south]{$\bar\sigma_z$};
        \filldraw[opacity=0.2, fill=black!20]
            plot[smooth, samples=50, domain=0:360] ( {2*cos(\x)}, {2*sin(\x)} );
        \draw[->] (0,0) -- (1.9696, 0.3473);
        \draw[->] (2.0285, -0.5435) arc (-15:0:2.1);
        \node[anchor=south] at (1, 0.15) {$r$};
        \node[anchor=west] at (2.0, -0.3) {$\beta$};
    \end{scope}
\end{tikzpicture}
    \caption{The three coordinates $(z, r, \beta)$ in the Lode coordinate system.}
    \label{fig:Lode}
\end{figure}

Moving from two-dimensional plane stress to full three-dimensional stress states increases the number of independent stress components from 3 to 6, which fundamentally raises the complexity of the problem. In the isotropic setting, the response can still be expressed in terms of principal stresses as illustrated in \Cref{fig:YieldSurface}, where sampling the three-dimensional yield surface in the Haigh--Westergaard space requires two angular coordinates $(\theta,\varphi)$ rather than the single angle $\theta$. There are also symmetry planes that simplify angular sampling. For instance, one may work with stress invariants instead of principal stresses and adopt the Lode coordinate system~\cite{Lode1926}, which is a cylindrical system with the three coordinates $(z, r, \beta)$. In this system, as shown in \Cref{fig:Lode}, the $z$-coordinate aligns with the hydrostatic line $\bar{\sigma}_x = \bar{\sigma}_y = \bar{\sigma}_z$, the radial coordinate $r$ lies in the deviatoric plane ($\pi$-plane), which is perpendicular to the $z$-axis, and the Lode angle $\beta$ describes the rotation around the $z$-axis within the deviatoric plane. Doing so, and having the freedom to assign any permutation of principal stresses to the three Cartesian axes in Haigh--Westergaard space, any arbitrary yield surface can be sampled at most within each $120^\circ$ sector of the Lode angle $\beta$ on a deviatoric plane. Consequently, the sampling domain has $181 \times 121 \times 50 \times 2 \approx 2.2\cdot10^{6}$ data points. However, at the same time, the number of Gauss points per element increases from 4 to 8 in three dimensions, so that the breakeven threshold~\cite{woldseth_use_2022} for a single increment with one Newton--Raphson iteration corresponds to a mesh of ${\sim}275,000$ elements (${\sim}65\times 65\times 65$ grid). This scenario, therefore, remains computationally favorable for a single simulation with a single increment and iteration.

The situation changes fundamentally when the count of stress components increases due to inelastic anisotropy. In this case, the constitutive response depends on the orientation of the stress tensor relative to the material axes, and principal stress states are no longer sufficient to characterize the material behavior. Considering the worst-case scenario with no symmetry at all, all six independent stress components must be accounted for, rendering grid-based sampling infeasible. This moves the breakeven threshold dramatically, and the surrogate model becomes economical only when a large number of simulations are to be performed.


The prohibitive data requirements of the fully anisotropic inelastic case can, however, be mitigated by more sophisticated sampling strategies. Latin hypercube sampling and related space-filling techniques offer a more efficient coverage of high-dimensional input spaces than structured grids. Physics-inspired sampling procedures~\cite{ghavamian_accelerating_2019,kalina_fetextrmann_2023} further reduce the sampling effort by tracking only the trajectories that are expected to be relevant for the intended application. Such targeted strategies come at the price of reduced robustness and applicability, as regions of the input space that were not sampled during training may lead to inaccurate predictions.

\section{Benchmark verification} \label{sec:Benchmark}
We verify the neural network surrogate by comparing its predictions against results from high-fidelity, concurrent FFT-based multiscale simulations. In this conventional approach, the macroscopic strain at each Gauss point is passed to a fine-scale RVE. Then the periodic boundary value problem governing the RVE, driven by this strain increment, is solved using the FFT-based homogenization scheme outlined in \Cref{alg:Procedure}. The resulting homogenized stress is then upscaled and returned to the macroscopic solver. This iterative two-scale process, which obviously incurs a prohibitive computational cost, is repeated for every Gauss point at each global iteration.

In contrast, the surrogate-equipped simulation bypasses this expensive scale transition entirely. For each Gauss point, a trial stress tensor, computed as the product of the strain increment and the homogenized stiffness tensor $\bar{\mathbb{C}}$, and the current equivalent plastic strain $\bar{\varepsilon}^p_{\text{eq}}$ are passed directly to the trained neural network. The network provides the updated stress and the incremental plastic multiplier $\Delta\bar{\gamma}$ in a single forward pass.

Three benchmark numerical examples were selected to comprehensively test the surrogate under varied stress states having different principal stress orientations. The examples are an internally pressurized annulus, an end-loaded cantilever beam, and Cook's membrane problem. To mitigate the extreme computational expense of the concurrent multiscale reference simulations, the macroscopic discretization was kept minimal. We employed eight-noded biquadratic elements with reduced integration, resulting in four Gauss points per element. This setup provides a tractable basis for comparing the accuracy and efficiency of the surrogate against the conventional multiscale approach.

We compare the CPU time of simulations using the neural network surrogate against that of the standard concurrent multiscale approach. Beyond raw computational speed, an even more compelling justification for using the surrogate as proposed in this paper lies in the severe memory constraints imposed by the conventional method. In our concurrent multiscale simulation, each Gauss point carries a staggering 1,048,576 internal variables ($512\times512$ pixels times 4 internal variables per pixel: three stress components and one plastic strain component), all of which must be stored and updated throughout the analysis. This places a severe restriction on the finite element discretization, as even on high-performance machines with hundreds of gigabytes of RAM, simulations with more than a handful of elements become infeasible. We perform all computations on a laptop running Arch Linux (kernel 6.18.2), equipped with an AMD Ryzen 7 PRO 6850U CPU operating at 2.7 GHz. For reference, generating the full training dataset with $1^\circ$ stress-angle sampling and the paired-increment strategy required 8 hours and 59 minutes, while training the network using PyTorch took only 18 seconds. Note that, as the conventional approach is far more computationally expensive, the comparison of simulation times is inherently one-sided if only the finite element analyses are considered. To provide a fair assessment, one must consider the total cost of the surrogate pipeline (data generation, network training, and finite element simulation) compared to that of the full concurrent multiscale modeling. This comparison demonstrates that, regardless of efficiency (see also~\Cref{sec:Scalability}), it is reasonable to train such networks even if one wants to use them for just a single analysis, as it enables fine discretizations that would otherwise be infeasible.

\subsection{Internally pressurized annulus}
\begin{figure}
    \centering
    \input{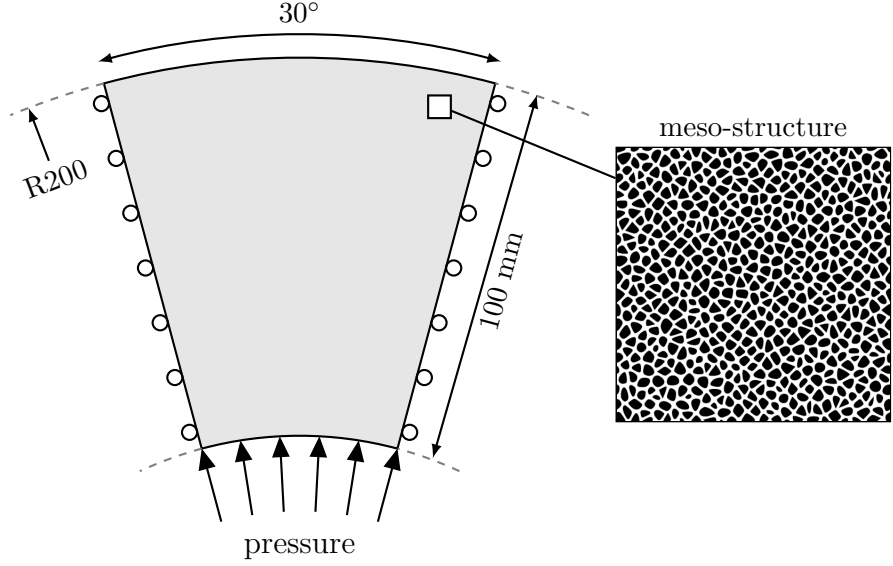}
    \caption{Schematic of the internally pressurized annulus problem and its meso-structure.}
    \label{fig:AnnulusSchematic}
\end{figure}

The first benchmark is an internally pressurized annulus with the geometry and boundary conditions illustrated in \Cref{fig:AnnulusSchematic}. The annulus has an inner radius of 100 mm, an outer radius of 200 mm, and is subjected to uniform internal pressure. Although the problem is axisymmetric, a $30^\circ$ segment is modeled with appropriate symmetric boundary conditions on the truncated edges. The segment is discretized into a 3$\times$3 mesh in the radial and circumferential directions, resulting in 9 elements in total. Pressure loading is applied indirectly within a displacement controlled regime by prescribing a total radial displacement of 4 mm at the inner edge.

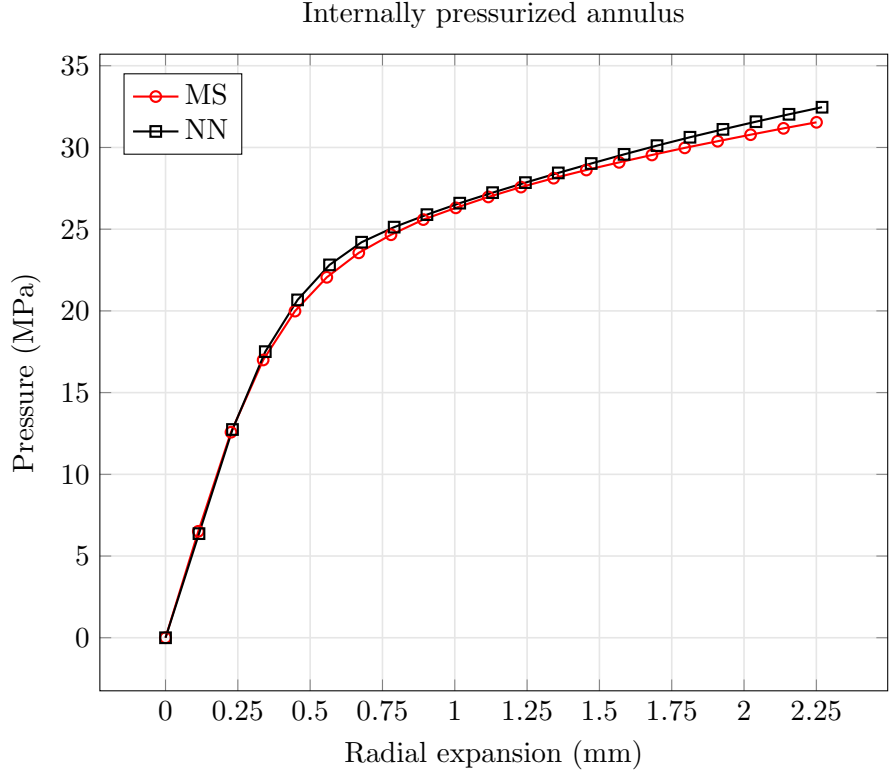
\begin{figure}
    \centering
    \begin{tikzpicture}
    \begin{axis}[
        width=12cm,
        height=10cm,
        grid=major,
        major grid style={line width=0.6pt, draw=black!10},
        xtick distance=0.25,
        xlabel={Radial expansion (mm)},
        ylabel={Pressure (MPa)},
        legend pos=north west,
        legend cell align=left,
        title=Internally pressurized annulus
    ]
        \addplot [red, thick, mark=o] table {data/mul-annulus.txt};
        \addlegendentry{MS}

        \addplot [black, thick, mark=square] table {data/dnn-annulus.txt};
        \addlegendentry{NN}
    \end{axis}
\end{tikzpicture}
    \caption{Pressure--radial expansion comparison between multiscale (MS) and neural network (NN) models.}
    \label{fig:AnnulusCurve}
\end{figure}

\begin{figure}
    \centering
    \input{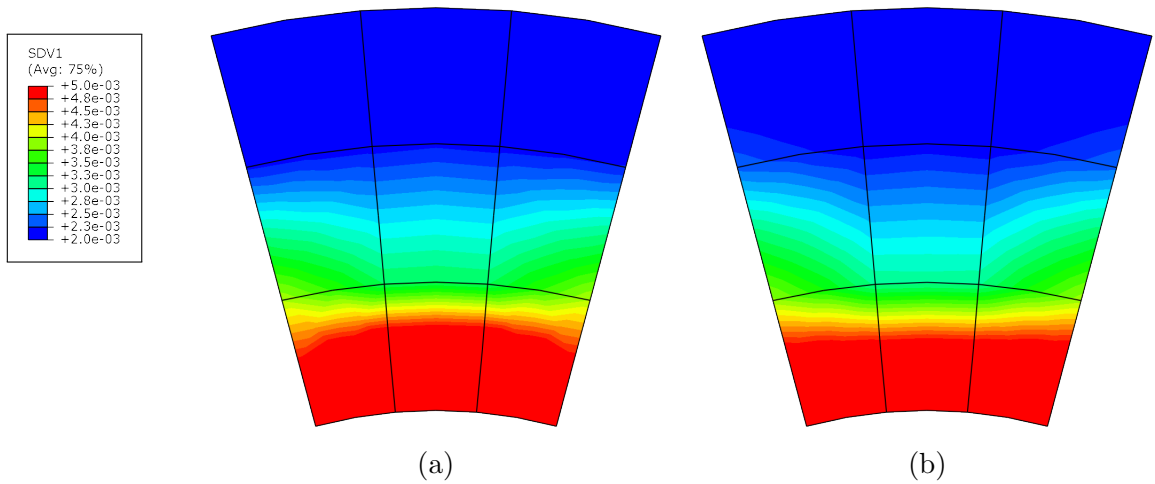}
    \caption{Plastic strain contours of annulus problem for (a) multiscale and (b) neural network models.}
    \label{fig:AnnulusContour}
\end{figure}

\Cref{fig:AnnulusCurve} compares the applied pressure versus the radial expansion at the outer edge for both the concurrent multiscale model and the neural network surrogate. \Cref{fig:AnnulusContour} shows the corresponding equivalent plastic strain distributions upon completion of the analyses. The results exhibit very good agreement. Minor deviations can be attributed to the subtle inherent anisotropy of the RVE in the high-fidelity concurrent multiscale simulation, which is approximated as fully isotropic in the surrogate model. Regarding computational time, the concurrent multiscale simulation required approximately 6 hours and 48 minutes, whereas the neural network surrogate completed the same analysis in approximately 0.8 seconds. This dramatic difference highlights the profound efficiency gain achieved by adopting the data-driven surrogate.

\subsection{End-loaded cantilever beam}
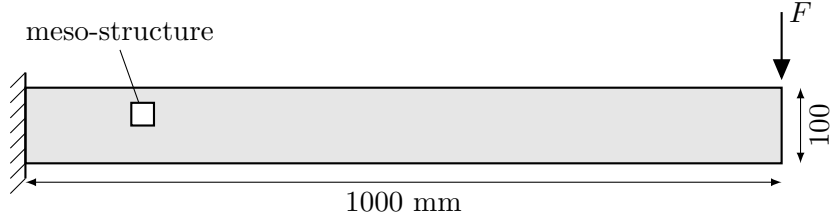
\begin{figure}
    \centering
    \begin{tikzpicture}
    \draw[fill=black!10, thick] (0,0) rectangle (10,1);
    \draw[->, >=triangle 45, thick, black] (10,2) -- (10,1.1);
    \node at (10.25,2) {$F$};

    \draw[thick] (0, -0.2) -- (0,1.2);
    \foreach \y in {-0.2, 0, ..., 1.2} {
        \draw[thin] (0,\y) -- (-0.2,\y-0.2);
    }

    \draw[<->, >=latex, black] (0,-0.25) -- (10,-0.25);
    \node at (5,-0.5) {1000 mm};
    \draw[<->, >=latex, black] (10.25,0) -- (10.25,1);
    \node[rotate=90] at (10.5,0.5) {100};
    
    \draw[fill=white, thick] (1.4,0.5) rectangle (1.7,0.8);
    \draw (1.5,0.8) -- (1.25,1.5);
    \node at (1.25,1.75) {meso-structure};
\end{tikzpicture}
    \caption{Schematic of the cantilever beam problem.}
    \label{fig:CantileverSchematic}
\end{figure}

The second benchmark is a cantilever beam subjected to a concentrated perpendicular load at its free end. The geometry and boundary conditions are illustrated in \Cref{fig:CantileverSchematic}. In contrast to the previous example, where the principal stress states were almost identical everywhere but the principal orientations varied spatially, this configuration generates a broad spectrum of stress states. Consequently, the stress angle $\theta$ varies from pure tension at the top edge, through pure shear at the neutral axis, to pure compression at the bottom edge. The beam is discretized with 2 elements through its thickness and 10 elements along its length, resulting in 20 elements in total.

To reduce computational cost, the first 10 elements (half of the beam) near the free end are assumed to remain elastic. Consequently, their response is computed by directly multiplying the strain with the homogenized stiffness tensor $\bar{\mathbb{C}}$. Only the remaining 10 elements near the fixed support follow the full concurrent multiscale modeling. We validated this simplification by a preliminary analysis confirming that the plastic zone does not extend into the elastic region. The problem is solved under displacement control by prescribing a vertical displacement of 200 mm at the neutral axis of the free end.

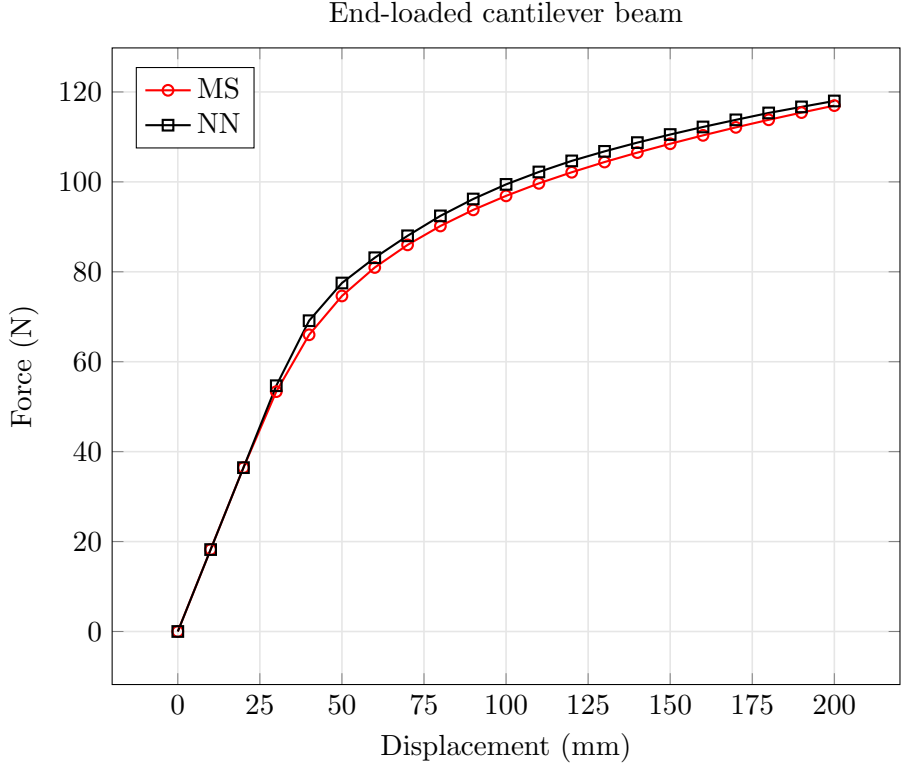
\begin{figure}
    \centering
    \begin{tikzpicture}
    \begin{axis}[
        width=12cm,
        height=10cm,
        grid=major,
        major grid style={line width=0.6pt, draw=black!10},
        xtick distance=25,
        xlabel={Displacement (mm)},
        ylabel={Force (N)},
        legend pos=north west,
        legend cell align=left,
        title=End-loaded cantilever beam
    ]
        \addplot [red, thick, mark=o] table {data/mul-cantilever.txt};
        \addlegendentry{MS}

        \addplot [black, thick, mark=square] table {data/dnn-cantilever.txt};
        \addlegendentry{NN}
    \end{axis}
\end{tikzpicture}
    \caption{Force--displacement comparison between multiscale (MS) and neural network (NN) models.}
    \label{fig:CantileverCurve}
\end{figure}

\begin{figure}
    \centering
    \input{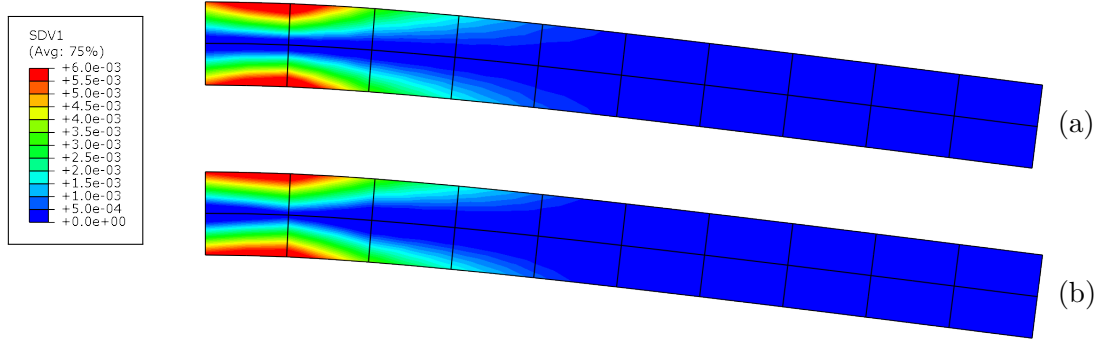}
    \caption{Plastic strain contours of cantilever beam problem for (a) multiscale and (b) neural network models.}
    \label{fig:CantileverContour}
\end{figure}

\Cref{fig:CantileverCurve} compares the resulting force--displacement curves for the concurrent multiscale model and the neural network surrogate. \Cref{fig:CantileverContour} presents the corresponding equivalent plastic strain contours. The curves show very good agreement, with minor deviations again attributable to the assumed isotropy in the surrogate model versus the subtle anisotropy in the high-fidelity RVE. The strain contours produced by both models are nearly identical, which further confirms the accuracy of the surrogate. Comparing the computation times reveals that the concurrent multiscale simulation required 7 hours and 18 minutes to complete, whereas the neural network surrogate required only 1.1 seconds. This result again demonstrates the substantial efficiency gain achieved by the data-driven approach over the conventional concurrent modeling paradigm.

\subsection{Cook's membrane problem}
\begin{figure}
    \centering
    \input{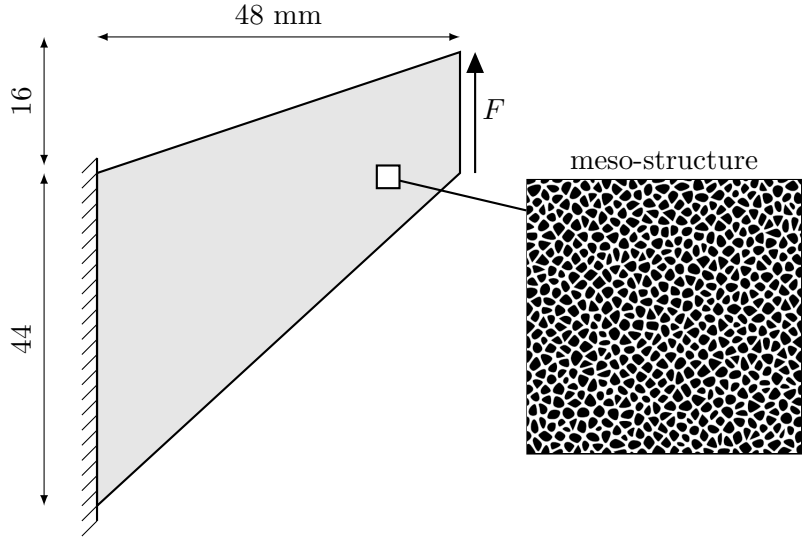}
    \caption{Schematic of Cook's problem and its meso-structure.}
    \label{fig:CooksSchematic}
\end{figure}

The third example is the classic Cook's membrane benchmark, widely used to evaluate finite element performance under combined bending and shear deformation~\cite{Neto2011}. As depicted in \Cref{fig:CooksSchematic}, the geometry is a tapered, skewed panel with unit thickness. The boundary conditions consist of a fully clamped left edge and a uniformly prescribed vertical displacement of 2 mm applied to the right edge. This example challenges the surrogate model with a diverse range of arbitrary stress states, spanning a broad spectrum of the stress angle $\theta$. The membrane is discretized using a 3$\times$3 mesh of eight-noded biquadratic elements, resulting in 9 elements in total.

\begin{figure}
    \centering
    \begin{tikzpicture}
    \begin{axis}[
        width=12cm,
        height=10cm,
        grid=major,
        major grid style={line width=0.6pt, draw=black!10},
        xtick distance=0.25,
        xlabel={Displacement (mm)},
        ylabel={Force (N)},
        legend pos=north west,
        legend cell align=left,
        title=Cook's membrane problem
    ]
        \addplot [red, thick, mark=o] table {data/mul-cooks.txt};
        \addlegendentry{MS}

        \addplot [black, thick, mark=square] table {data/dnn-cooks.txt};
        \addlegendentry{NN}
    \end{axis}
\end{tikzpicture}
    \caption{Force--displacement comparison between multiscale (MS) and neural network (NN) models.}
    \label{fig:CooksCurve}
\end{figure}
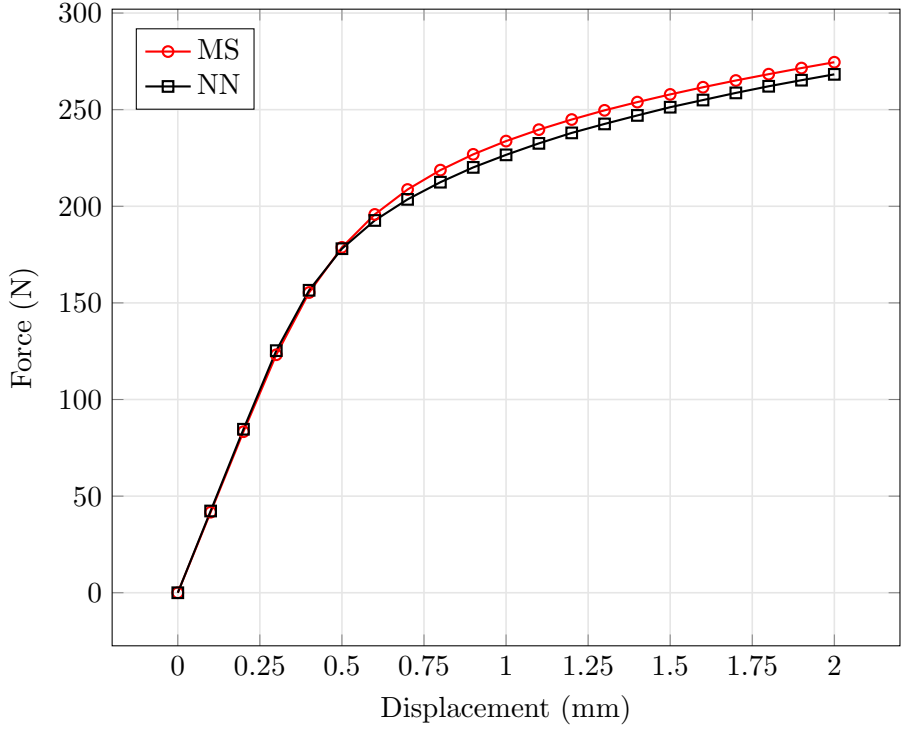

\begin{figure}
    \centering
    \input{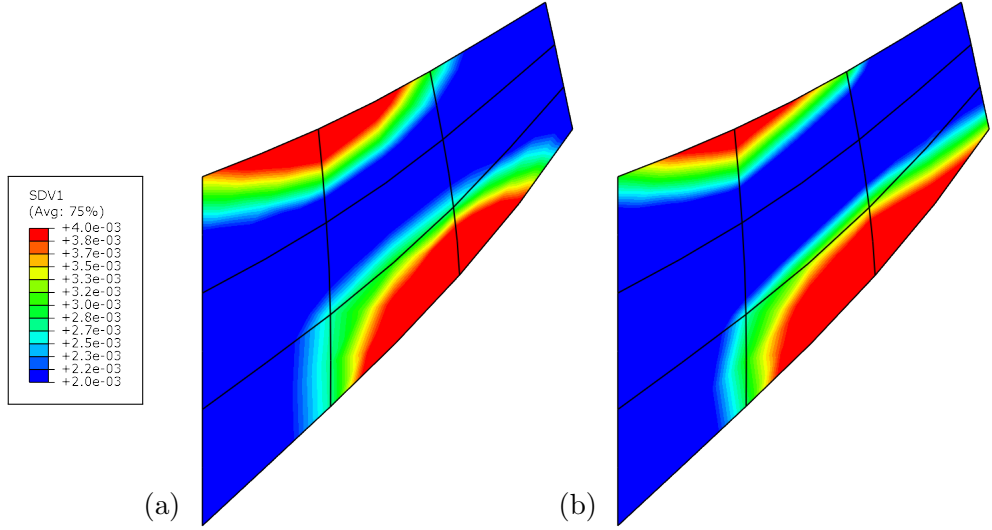}
    \caption{Plastic strain contours of Cook's problem for (a) multiscale and (b) neural network models.}
    \label{fig:CooksContour}
\end{figure}

The resulting force--displacement curves for both the multiscale model and the neural network surrogate are compared in \Cref{fig:CooksCurve}. The displacement is measured at the upper right vertex of the membrane. The curves again exhibit very good agreement, with minor deviations stemming from the isotropic approximation in the surrogate versus the slight anisotropy in the high-fidelity simulation. The equivalent plastic strain contours, presented for both models in \Cref{fig:CooksContour}, are nearly identical. In terms of computational efficiency, the concurrent multiscale simulation required 5 hours and 13 minutes, while the neural network surrogate completed the same analysis in 0.8 seconds. This dramatic reduction in computational time highlights both the practical advantage of the proposed data-driven approach and the necessity of the paradigm shift from conventional methods to surrogate-based constitutive modeling.

\section{Sensitivity analysis} \label{sec:sensitivity}
This section presents three groups of sensitivity analyses to assess the robustness and data requirements of the surrogate model. In the first group, the network trained on the complete dataset is evaluated under varying numerical conditions, including mesh size and loading step size. The mesh sensitivity analysis examines the effect of mesh refinement, which increases the diversity of local stress states and provides a more comprehensive sampling of the stress angle $\theta$ spectrum. The step-size sensitivity analysis, on the other hand, examines whether the paired-increment strategy successfully enables the network to handle strain increments across the full range spanned by the training pairs.

The second group investigates the sensitivity of the network to the size and composition of the training data. For this purpose, separate networks are trained on a reduced dataset, and their predictions are compared against those of the model trained on the complete dataset. It is important to note that the same neural network architecture employed for the model trained on the full dataset in the previous tests is used here as well. First, the density of stress-angle sampling is varied by generating datasets with coarser increments of $2^\circ$ and $5^\circ$ in addition to the baseline $1^\circ$. This allows us to quantify how much the data generation effort can be reduced before compromising prediction accuracy. We then isolate the role of the paired-increment strategy by training a network exclusively on the larger increments and omitting the initial small steps. This comparison directly assesses whether the paired strategy is essential for equipping the network with the robustness to deal with a continuous spectrum of strain increment magnitudes. 

Finally, in the third group, the model's ability to generalize to out-of-distribution loading is evaluated using loading paths that were not included in the training dataset. These cases include non-monotonic loading histories involving repeated loading and unloading, as well as non-proportional loading paths where the loading direction changes during the deformation process.

\begin{figure}
    \centering
    \input{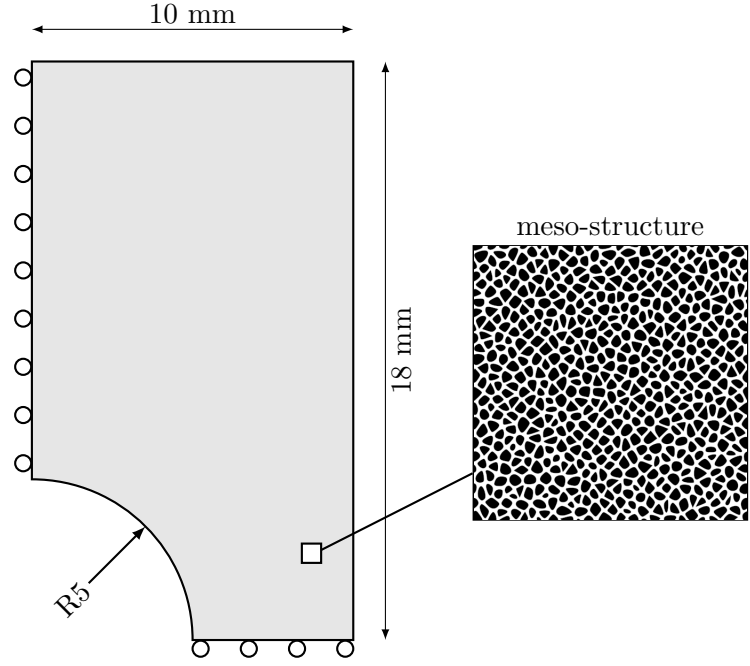}
    \caption{Schematic of the perforated plate problem and its meso-structure.}
    \label{fig:PerforatedSchematic}
\end{figure}

We employ the classical perforated plate problem, a common benchmark for plane-stress elasto-plasticity. The plate dimensions are $20 \times 36$ mm, with a central circular hole of 10 mm diameter. For all cases except the non-proportional loading, the plate is loaded in uniaxial tension via a uniform prescribed displacement applied to the upper and lower edges. For the latter, to obtain a non-proportional loading path, we first apply a uniform pressure over the face of the hole, and then impose a uniformly distributed tensile load over the top edge of the plate in a subsequent step. Due to symmetry in both geometry and loading, only one quarter of the plate is modeled, as illustrated in \Cref{fig:PerforatedSchematic}. Appropriate symmetric boundary conditions are enforced on the truncated edges. 

\subsection{Mesh sensitivity}
\begin{figure}
    \centering
    \begin{tikzpicture}
    \begin{axis}[
        width=12cm,
        height=10cm,
        grid=major,
        major grid style={line width=0.6pt, draw=black!10},
        xtick distance=0.05,
        xlabel={Displacement (mm)},
        ylabel={Force (MPa)},
        legend pos=north west,
        legend cell align=left,
        title=Mesh sensitivity,
        tick label style={/pgf/number format/fixed}
    ]
        \addplot [red, thick, mark=none] table {data/perforated-mesh1.txt};
        \addlegendentry{Mesh 1}

        \addplot [black, thick, mark=none, dashed] table {data/perforated-mesh2.txt};
        \addlegendentry{Mesh 2}

        \addplot [black, thick, mark=square, only marks] table {data/perforated-mesh3.txt};
        \addlegendentry{Mesh 3}
    \end{axis}
\end{tikzpicture}
    \caption{Force--displacement response for the three meshes in the mesh sensitivity analysis.}
    \label{fig:MeshCurve}
\end{figure}
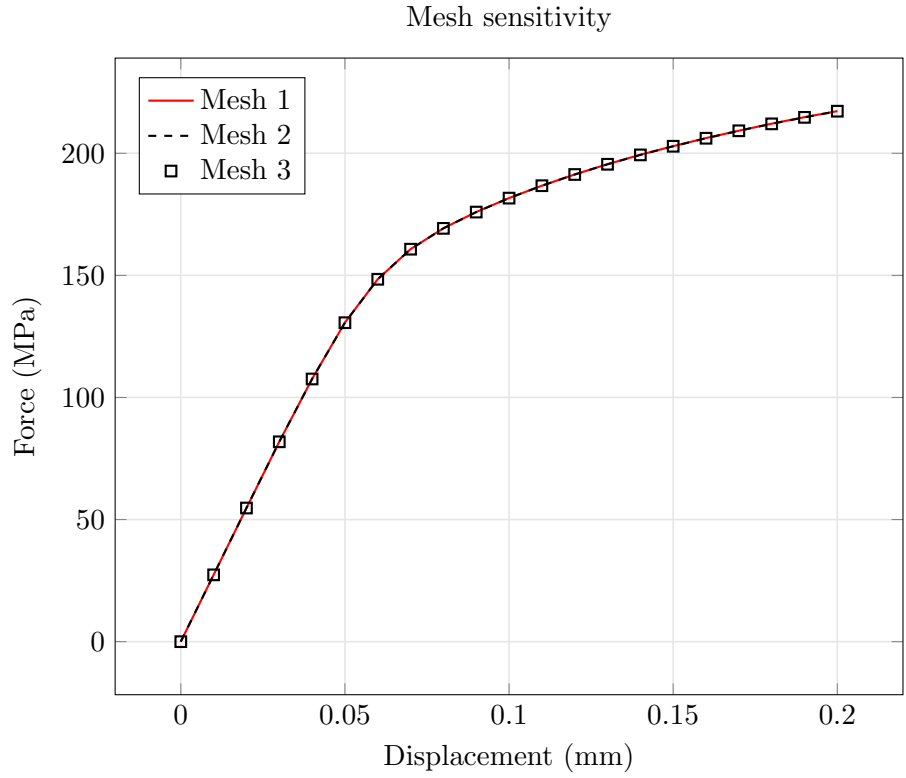

\begin{figure}
    \centering
    \input{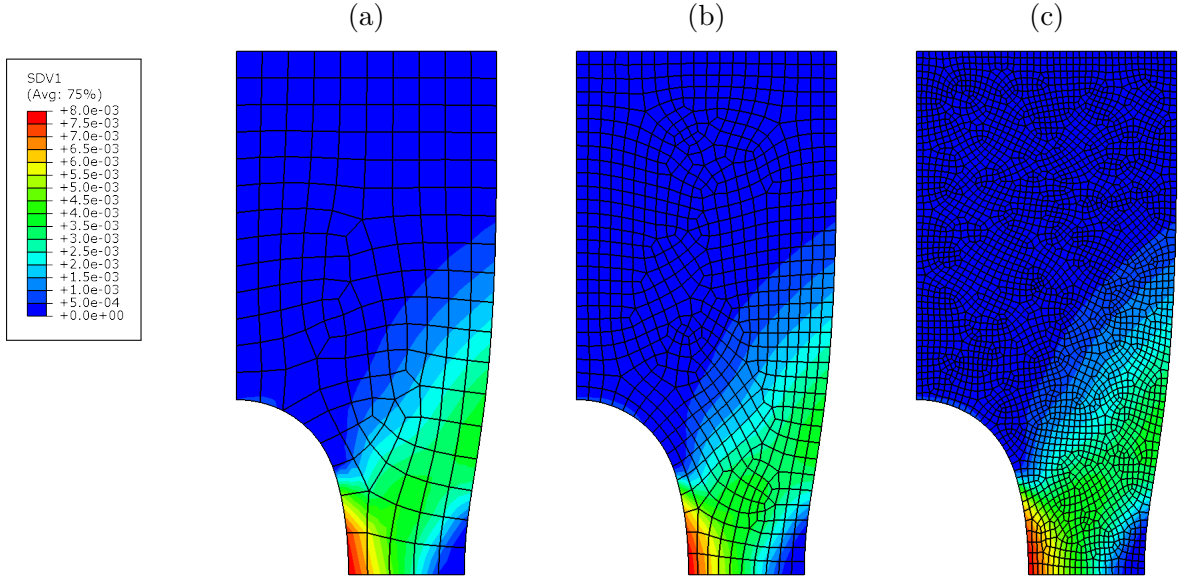}
    \caption{Equivalent plastic strain contours for the mesh sensitivity analysis: (a) Mesh 1, (b) Mesh 2, and (c) Mesh 3.}
    \label{fig:SensitivityContour}
\end{figure}

For the mesh sensitivity analysis, three discretizations, including a coarse, a medium, and a fine mesh using reduced-integration eight-noded biquadratic elements, are considered. The coarse mesh is generated with a seed spacing of 1 mm, resulting in elements with an approximate side length of 1 mm. The medium mesh uses half that spacing (0.5 mm), and the fine mesh a quarter of the original spacing (0.25 mm). Consequently, the meshes, designated as Mesh 1, Mesh 2, and Mesh 3, comprise 177 elements, 760 elements, and 3146 elements, respectively. \Cref{fig:MeshCurve} presents a comparison of the total applied force over half of one edge versus the prescribed displacement for all three meshes. The curves are virtually indistinguishable. The equivalent plastic strain contour for each mesh is also presented in \Cref{fig:SensitivityContour}, showing a similar degree of agreement. It must be mentioned here that, from this point onward, the medium mesh, generated by 760 eight-noded biquadratic elements, is the basis for the remaining sensitivity analyses.

\subsection{Step-size sensitivity}
In nonlinear finite element analysis, automatic load stepping is the standard practice to ensure solution robustness and efficiency. This algorithm adapts the load increment size based on the number of Newton--Raphson iterations required for convergence. If the iteration count exceeds a specified limit, the increment is automatically subdivided, and if convergence is achieved too rapidly, the increment size is increased. Commercial solvers like Abaqus implement this capability as their default strategy for incremental nonlinear analysis.

In this part of the sensitivity analysis, we assess the robustness of the neural network surrogate with respect to increment size, i.e., its ability to provide consistent stress updates across a range of strain increments without solution deviation. To perform this step-size sensitivity analysis, we temporarily disable two key adaptive features. First, we suppress the UMAT signal that would request increment halving based on the incremental plastic multiplier, and second, we switch from automatic to fixed load increments. The solver is also configured to permit more iterations than usual to avoid premature termination. Doing so, we isolate the effect of increment size.

For the perforated plate benchmark, we target a maximum vertical displacement of 0.2 mm at the upper edge and apply it using three different fixed increments, namely 50, 25, and 10 steps. The simulation with 10 steps diverged, so we conducted a manual search for a stable step count. The smallest number of fixed steps that converged was 16. It is important to note that re-enabling adaptive stepping, along with allowing the UMAT to signal for increment subdivision, would prevent such divergence, but we intentionally disable these features here to evaluate the performance of the neural network under a purely fixed scheme. However, this non-adaptive approach is not recommended for general nonlinear analysis.

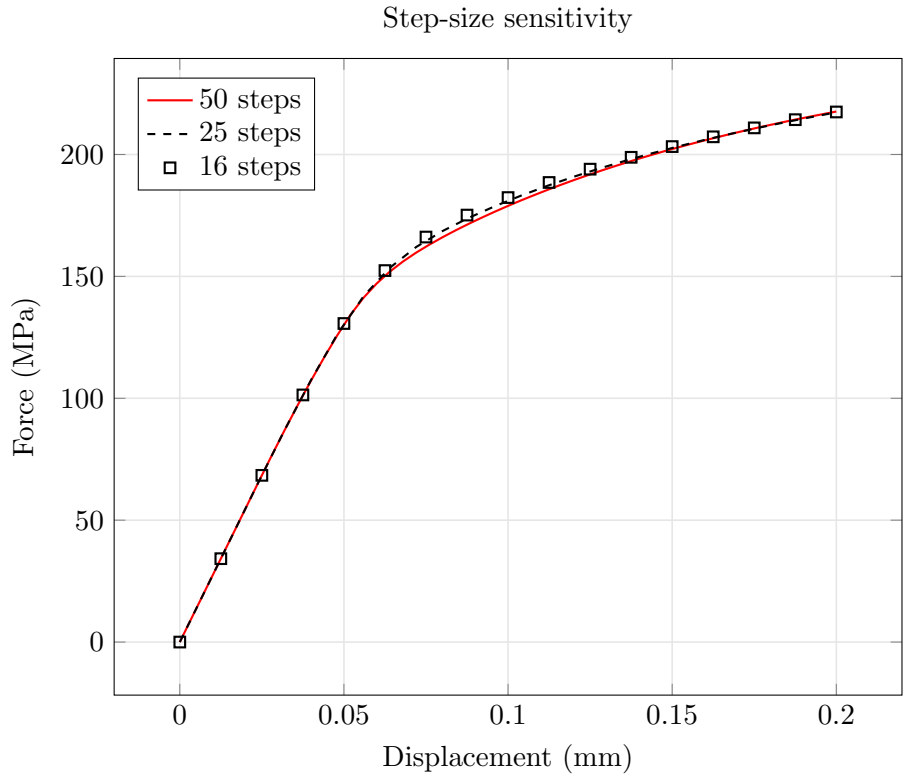
\begin{figure}
    \centering
    \begin{tikzpicture}
    \begin{axis}[
        width=12cm,
        height=10cm,
        grid=major,
        major grid style={line width=0.6pt, draw=black!10},
        xtick distance=0.05,
        xlabel={Displacement (mm)},
        ylabel={Force (MPa)},
        legend pos=north west,
        legend cell align=left,
        title=Step-size sensitivity,
        tick label style={/pgf/number format/fixed}
    ]
        \addplot [red, thick, mark=none] table {data/perforated-step1.txt};
        \addlegendentry{50 steps}

        \addplot [black, thick, mark=none, dashed] table {data/perforated-step2.txt};
        \addlegendentry{25 steps}

        \addplot [black, thick, mark=square, only marks] table {data/perforated-step3.txt};
        \addlegendentry{16 steps}
    \end{axis}
\end{tikzpicture}
    \caption{Force--displacement response for the three converged step sizes in the step-size sensitivity analysis.}
    \label{fig:StepCurve}
\end{figure}

The resulting force--displacement curves for the converged cases (50, 25, and 16 steps) are compared in \Cref{fig:StepCurve}. The curves show very good agreement, indicating that the trained neural network provides consistent and accurate stress updates across a practical range of fixed increment sizes.

\subsection{Stress-angle sampling}
\begin{figure}
    \centering
    \begin{tikzpicture}
    \begin{axis}[
        width=12cm,
        height=10cm,
        grid=major,
        major grid style={line width=0.6pt, draw=black!10},
        xtick distance=0.05,
        xlabel={Displacement (mm)},
        ylabel={Force (MPa)},
        legend pos=north west,
        legend cell align=left,
        title=Stress-angle sampling,
        tick label style={/pgf/number format/fixed}
    ]
        \addplot [red, thick, mark=none] table {data/perforated-theta1.txt};
        \addlegendentry{1$^\circ$ increments}

        \addplot [black, thick, mark=none] table {data/perforated-theta2.txt};
        \addlegendentry{2$^\circ$ increments}

        \addplot [black, thick, mark=none, dashed] table {data/perforated-theta3.txt};
        \addlegendentry{5$^\circ$ increments}
    \end{axis}
\end{tikzpicture}
    \caption{Comparison of force--displacement responses for training datasets with different stress-angle resolutions.}
    \label{fig:ThetaCurve}
\end{figure}
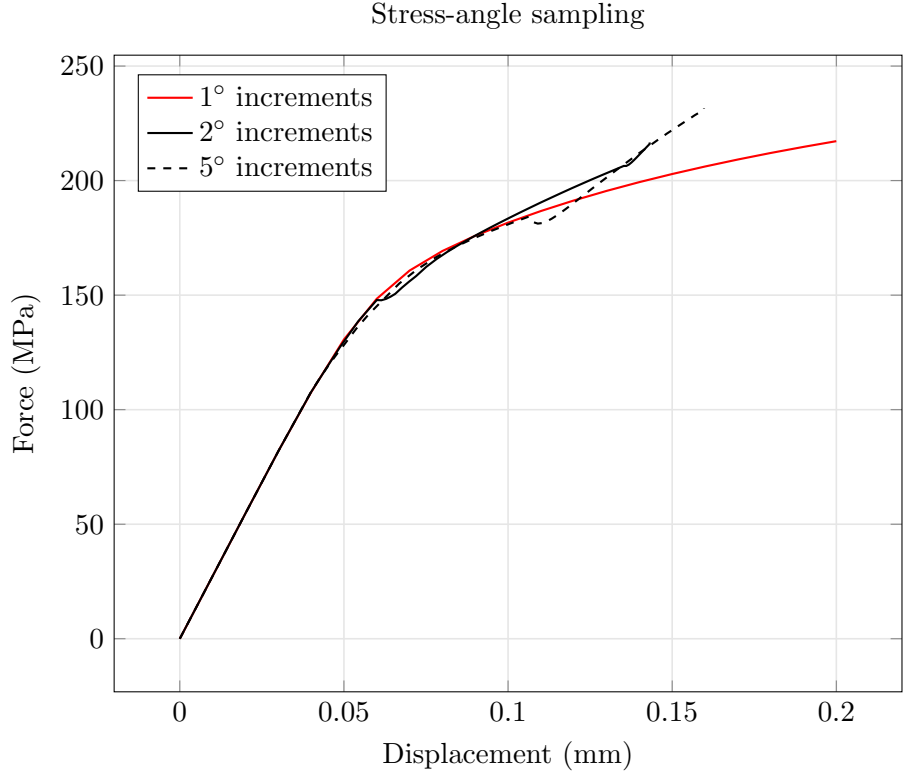

This section investigates the sensitivity of the surrogate model to the density of stress-angle sampling within the training data. We compare three sampling strategies, defined by stress-angle increments of $1^\circ$, $2^\circ$, and $5^\circ$. The first one, which constitutes the full dataset and serves as the baseline, was used for all previous analyses. For the coarser $2^\circ$ and $5^\circ$ increments, separate neural networks are trained on the correspondingly reduced datasets. The resulting force--displacement curves are presented in \Cref{fig:ThetaCurve}. As evident from the figure, the response curves produced by the models trained on the reduced datasets deviate from the baseline analysis. Moreover, these coarser sampling strategies produce erroneous stress updates, which lead to premature simulation divergence. These results indicate that a $1^\circ$ sampling increment appears to be necessary for ensuring accuracy, and it serves as an effective practical threshold for robust performance across all stress angles.

\subsection{Paired-increment strategy}
\begin{figure}
    \centering
    \begin{tikzpicture}
    \begin{axis}[
        width=12cm,
        height=10cm,
        grid=major,
        major grid style={line width=0.6pt, draw=black!10},
        xtick distance=0.05,
        xlabel={Displacement (mm)},
        ylabel={Force (MPa)},
        legend pos=north west,
        legend cell align=left,
        title=Paired-increment strategy,
        tick label style={/pgf/number format/fixed}
    ]
        \addplot [red, thick, mark=none] table {data/perforated-paired1.txt};
        \addlegendentry{Reference}

        \addplot [black, thick, mark=none] table {data/perforated-paired2.txt};
        \addlegendentry{50 steps}

        \addplot [black, thick, mark=none, dashed] table {data/perforated-paired3.txt};
        \addlegendentry{25 steps}

        \addplot [black, thick, mark=square, only marks] table {data/perforated-paired4.txt};
        \addlegendentry{16 steps}
    \end{axis}
\end{tikzpicture}
    \caption{Comparison of force--displacement responses from the surrogate model trained on the reduced dataset (different step counts) against the reference curve from the model using the full dataset.}
    \label{fig:PairedCurve}
\end{figure}
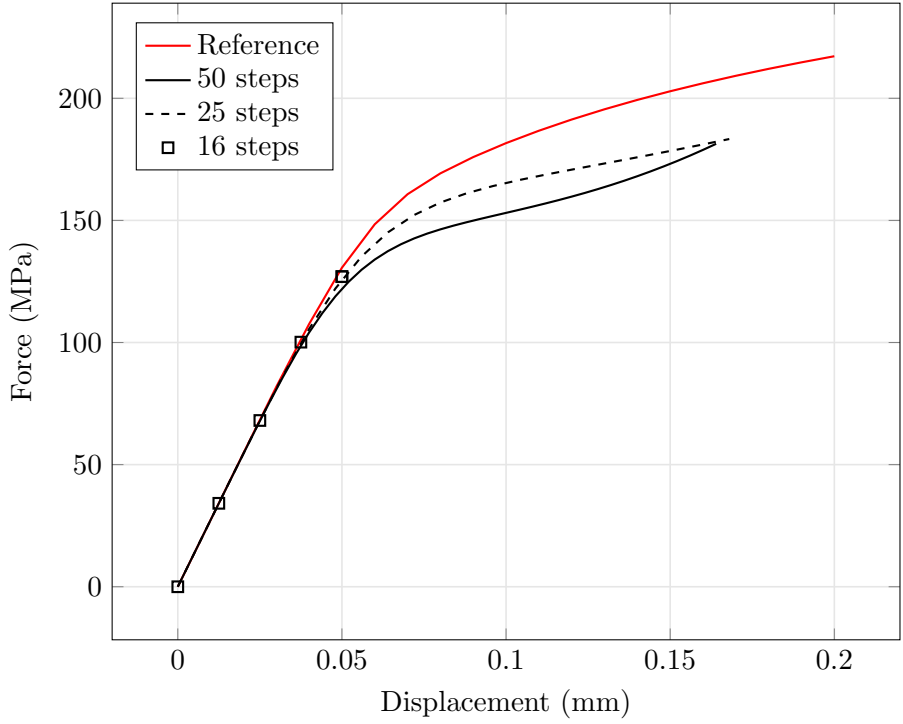

This sensitivity analysis evaluates the effectiveness of the proposed paired-increment strategy. To this end, a separate neural network is trained on a reduced dataset that retains only the larger plastic increment and omits the initial small elastic-violating step from each pair. The trained network is then tested under the same three fixed loading schemes (50, 25, and 16 steps) used in the step-size sensitivity analysis. The resulting force--displacement curves are compared against a reference solution, obtained using the network trained on the full paired-increment dataset, in \Cref{fig:PairedCurve}.

Simulations using the reduced-strategy network failed to complete the loading path, as all three cases diverged at an intermediate step. Moreover, even before divergence, the predicted responses exhibited significant deviation from the reference solution. This outcome demonstrates that the paired-increment strategy, which includes both the initial small step that barely violates the elastic limit and the subsequent larger step that induces noticeable plastic deformation, is essential. This strategy provides the network with the necessary information to accurately perform stress updates under arbitrary magnitudes of strain increments across a wide, yet numerically reasonable, range.

\subsection{Out-of-distribution loading}
\subsubsection{Non-monotonic}
\begin{figure}
    \centering
    \begin{tikzpicture}
    \begin{axis}[
        width=12cm,
        height=10cm,
        grid=major,
        major grid style={line width=0.6pt, draw=black!10},
        xtick distance=0.05,
        xlabel={Displacement (mm)},
        ylabel={Force (MPa)},
        legend pos=north west,
        legend cell align=left,
        title=Non-monotonic loading,
        tick label style={/pgf/number format/fixed}
    ]
        \addplot [red, thick, mark=none] table {data/perforated-mono.txt};
        \addlegendentry{Monotonic}
        
        \addlegendimage{only marks, mark=square, thick, black}
        \addlegendentry{Case 1}
        \addlegendimage{only marks, mark=o, thick, black}
        \addlegendentry{Case 2}
        \addlegendimage{only marks, mark=triangle, mark size=3pt, thick, black}
        \addlegendentry{Case 3}

        \addplot [black, thick, loosely dashed, mark=square, mark options={solid, draw=black, fill=none}, mark repeat=27, mark phase=26] table {data/perforated-cyclic1.txt};
        
        \addplot [black, thick, loosely dashed, mark=o, mark options={solid, draw=black, fill=none}, mark repeat=36, mark phase=35] table {data/perforated-cyclic2.txt};
        
        \addplot [black, thick, loosely dashed, mark=triangle, mark size=3pt, mark options={solid, draw=black, fill=none}, mark repeat=36, mark phase=35] table {data/perforated-cyclic3.txt};
    \end{axis}
\end{tikzpicture}
    \caption{Comparison of force--displacement responses for monotonic and non-monotonic loading cases.}
    \label{fig:NonMonotonicCurve}
\end{figure}
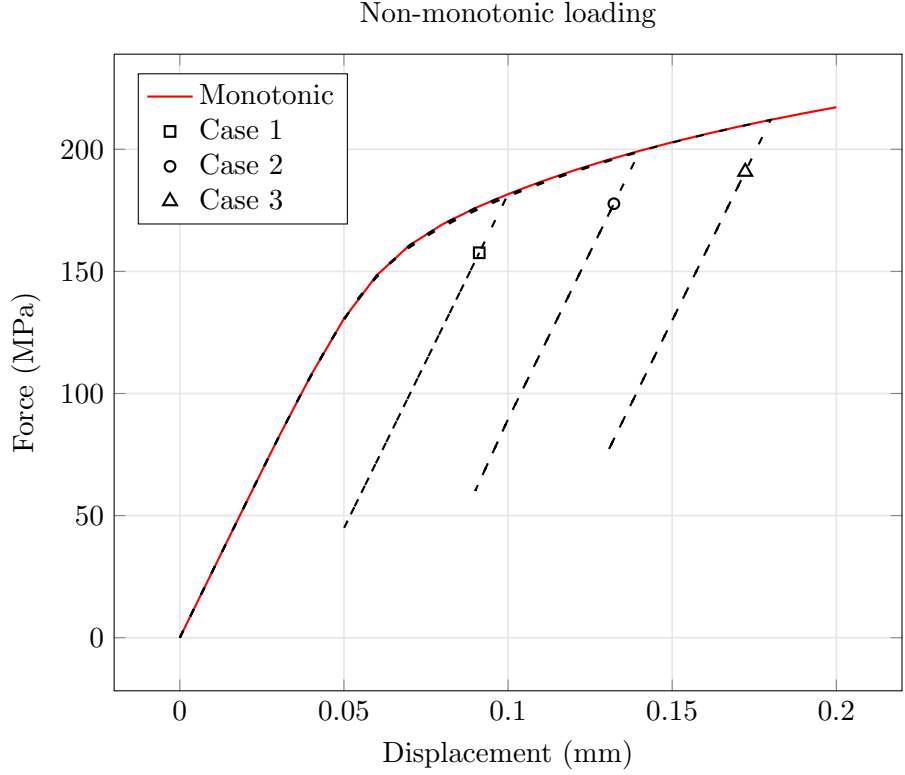

For the non-monotonic loading case, repeated loading and unloading are applied. A prescribed displacement of 0.1 mm is first imposed on the upper edge of the plate, after which the plate is unloaded to half of this displacement and subsequently reloaded beyond the initial 0.1 mm. The model accurately predicts both the initial loading and the subsequent unloading. However, during reloading, it fails to converge at the transition from elastic to elastoplastic response. The same procedure is repeated with initial prescribed displacements of 0.14 mm and 0.18 mm, yielding identical behavior. The resulting force--displacement curves for all three cases are compared with the corresponding monotonic loading response in \Cref{fig:NonMonotonicCurve}. The three cases in the figure correspond to initial prescribed displacements of 0.10, 0.14, and 0.18 mm. The points at which the analyses diverge are indicated by square, circular, and triangular markers, respectively.

These results indicate that, although the model can accurately predict unloading, it is unable to represent the stress update correctly during reloading. This behavior stems from the formulation adopted in the present work, where the sign of the plastic strain increment is used as an indicator of plastic evolution. When the plastic strain increment is positive, the updated stress state provided by the network is the valid one, whereas for a negative increment the trial state is retained, which naturally captures unloading. Reloading, however, requires the model to identify the transition from elastic to elastoplastic behavior for material points with an existing plastic history. Since the training dataset only contains this transition for the initial onset of yielding and does not include corresponding transitions following prior plastic deformation, the network is unable to predict the stress evolution accurately in this regime. Consequently, erroneous stress updates lead to loss of convergence.

\subsubsection{Non-proportional}
For the non-proportional loading case, the same perforated plate geometry is employed, but two distinct loading stages are applied to generate different loading paths. In the first stage, a uniform pressure is applied to the inner surface of the hole. In the second stage, this pressure is maintained while a uniformly distributed tensile load is applied to the upper edge of the plate. Three loading cases are considered, with inner hole pressures of 16, 18, and 20 MPa, followed by the same uniformly distributed tensile load of 10 MPa in the second stage.

\begin{figure}
    \centering
    \input{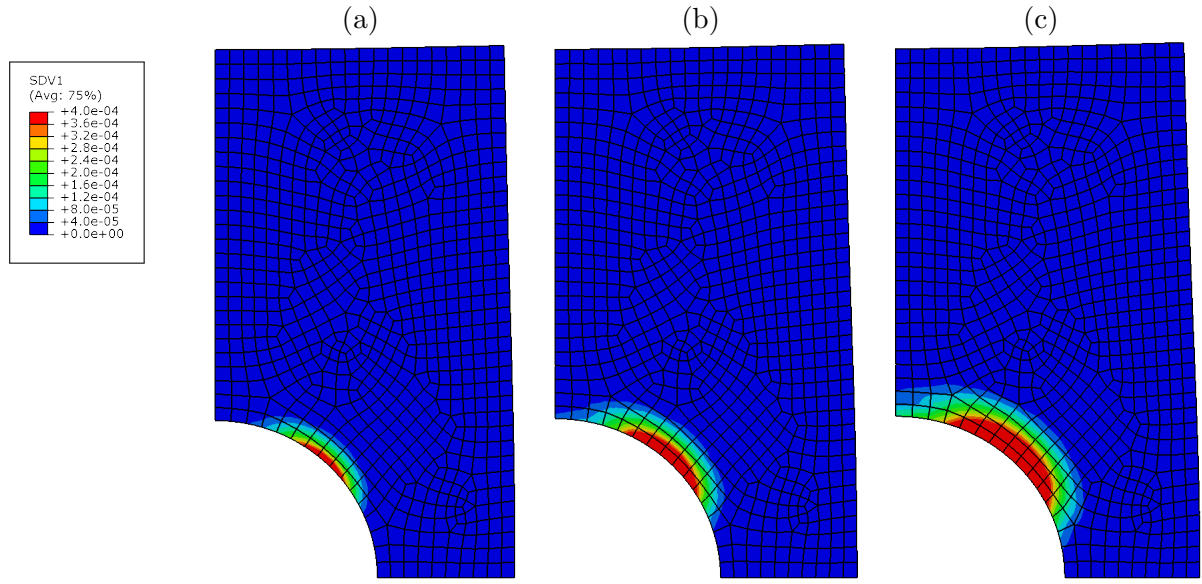}
    \caption{Equivalent plastic strain contours for the non-proportional loading case at the end of the first loading stage: (a) Case 1, (b) Case 2, and (c) Case 3.}
    \label{fig:NonProportionalContourA}
\end{figure}

The equivalent plastic strain contours at the end of the first loading stage are shown in \Cref{fig:NonProportionalContourA} to illustrate the extent of plastic deformation produced by the applied hole pressure. As shown, the 16 MPa case produces only a very small plastic zone. Increasing the pressure to 18 MPa results in a noticeable plastic region, while the 20 MPa case leads to substantial plastic deformation around the hole.

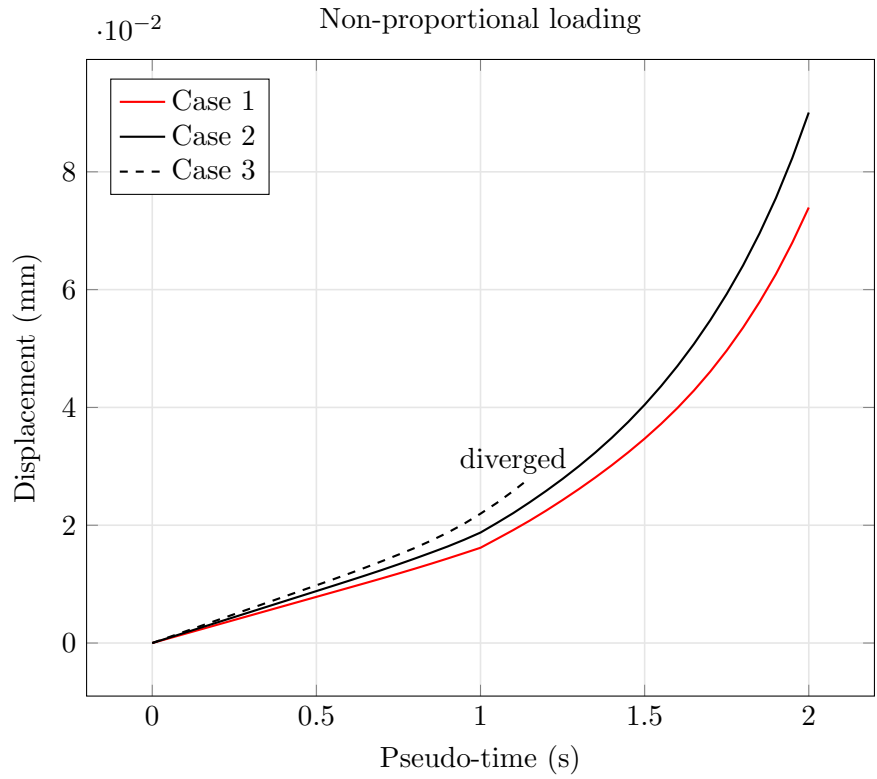
\begin{figure}
    \centering
    \begin{tikzpicture}
    \begin{axis}[
        width=12cm,
        height=10cm,
        grid=major,
        major grid style={line width=0.6pt, draw=black!10},
        xtick distance=0.5,
        xlabel={Pseudo-time (s)},
        ylabel={Displacement (mm)},
        legend pos=north west,
        legend cell align=left,
        title=Non-proportional loading,
        tick label style={/pgf/number format/fixed}
    ]
        \addplot [red, thick, mark=none] table {data/perforated-non1.txt};
        \addlegendentry{Case 1}

        \addplot [black, thick, mark=none] table {data/perforated-non2.txt};
        \addlegendentry{Case 2}

        \addplot [black, thick, mark=none, dashed] table {data/perforated-non3.txt};
        \addlegendentry{Case 3}

        \node[scale=1.0] (label) at (1.1, 0.031) {diverged};

    \end{axis}
\end{tikzpicture}
    \caption{Vertical displacement of the node located at the top of the hole ($t=0$ to $t=1$ corresponds to the first loading stage and $t=1$ to $t=2$ corresponds the second loading stage).}
    \label{fig:NonProportionalCurve}
\end{figure}

\begin{figure}
    \centering
    \input{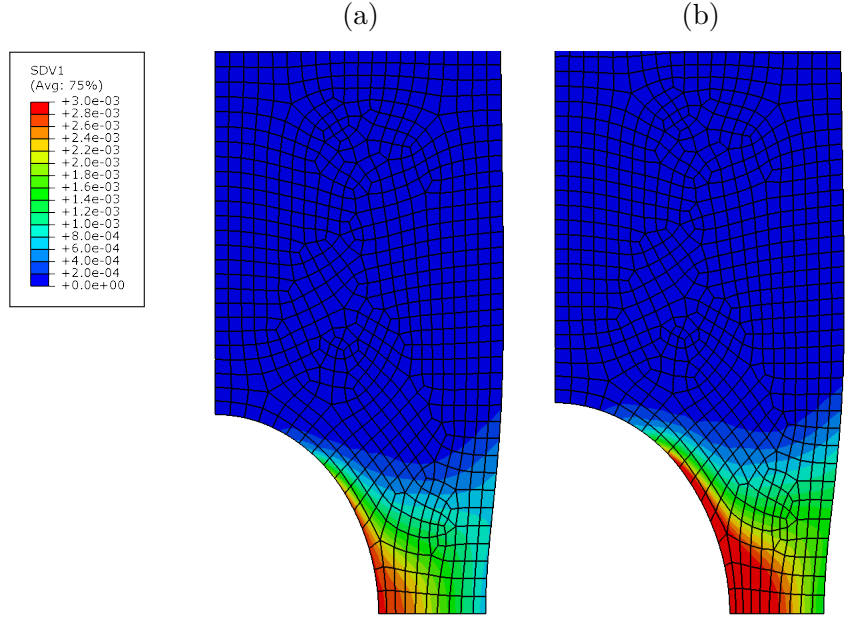}
    \caption{Equivalent plastic strain contours for the non-proportional loading case at the end of the analysis: (a) Case 1 and (b) Case 2.}
    \label{fig:NonProportionalContourB}
\end{figure}

During the second loading stage, the first two cases converge successfully, whereas the third case, corresponding to an inner hole pressure of 20 MPa, fails to converge shortly after the second stage begins. This behavior is illustrated in \Cref{fig:NonProportionalCurve}, where the vertical displacement of the node located at the top of the hole ($\theta=90^\circ$) is plotted against pseudo-time. Here, $t=0$ to $t=1$ corresponds to the first loading stage (hole pressure), while $t=1$ to $t=2$ represents the second loading stage (tensile loading of the upper edge). The equivalent plastic strain contours at the end of the analysis are also presented in \Cref{fig:NonProportionalContourB}.

The results indicate that the first two cases converge because the plastic regions generated during the two loading stages remain largely separated. Consequently, most material points that undergo plastic deformation during the second stage have not experienced noticeable prior plastic yielding, and the network encounters conditions not too far from those represented in the training dataset. In contrast, for the 20 MPa case, the plastic region generated during the first stage extends significantly to the right of the hole, causing a substantial overlap with the plastic zone produced during the second stage. As a result, material points with an existing plastic history continue to evolve plastically along a different loading path than those represented in the training data. Since such loading histories are absent from the training dataset, the network produces inaccurate stress updates, ultimately leading to loss of convergence during the early increments of the second loading stage. While such a failure to converge might seem inconvenient at first, it is actually a benefit regarding robustness: The setup does not simply continue with wrong results, but indicates when it gives wrong predictions.

\section{Concluding remarks} \label{sec:conclusion}
This work presented a deep learning surrogate model to overcome the main computational bottleneck of concurrent multiscale plasticity, stemming from the prohibitive cost of solving nested fine-scale boundary value problems at each Gauss point. The objective was to replace this expensive procedure with a data-driven constitutive update performed by a neural network that directly learns the constitutive response of a heterogeneous RVE from training data generated via FFT-based homogenization.

Three benchmark studies demonstrated that the neural network surrogate reproduces the responses of the high-fidelity concurrent multiscale simulations with very good agreement. Minor discrepancies were attributed to the subtle anisotropy inherent in the RVE in the reference simulations, which is approximated as isotropic in the surrogate formulation. Most notably, the computational time was reduced from several hours to approximately one second per analysis, achieving speed-ups of 23,000--30,000. These results clearly establish the surrogate model as a viable and highly efficient alternative to conventional concurrent multiscale plasticity.

Sensitivity studies verified perfect mesh objectivity and revealed only very minor sensitivity to the chosen step size, which is also frequently seen in conventional return-mapping approaches. The investigation of stress-angle sampling revealed that sufficiently dense angular coverage consisting of $1^\circ$ increments is essential for stable and accurate stress updates. In addition, the paired-increment strategy is crucial for enabling the network to generalize across a continuous spectrum of strain increments. Omitting this strategy leads to deviation from the reference results and premature divergence. The surrogate is trained on monotonic, proportional loading paths, so its predictions degrade for loading outside this distribution, e.g., non-monotonic or non-proportional paths. Within this scope, however, the intended benefit, replacing the fine-scale solve for monotonic-proportional loading, is fully realized. 

Herein also lies the principal limitation. The stress space must be sampled with sufficient density. For the two-dimensional isotropic case, a modest dataset of $18,100$ points suffices. Extending to three-dimensional isotropic plasticity increases this by a factor of ${\sim}120$. Yet the surrogate remains economical even for a single simulation, as the Gauss point evaluations in a three-dimensional mesh of reasonable size quickly exceed this number. For anisotropic plasticity, however, the dataset grows substantially, and the breakeven threshold shifts dramatically (see~\Cref{sec:Scalability} for details).  However, such an investment can be justified, as the resulting surrogate constitutes a general-purpose constitutive model applicable across a broad class of simulations.

Beyond raw computational speed, an even more compelling justification for the surrogate lies in the severe memory constraints imposed by the conventional method. In our simulations, each Gauss point carries over a million internal variables, rendering concurrent simulations even with very coarse meshes infeasible. This demonstrates that training such networks is worthwhile even for a single analysis, as it unlocks simulations that lie far beyond the reach of traditional concurrent multiscale frameworks.

\section{Declaration of Competing Interests} \label{sec:Declaration}
The authors declare no conflicts of interest.

\section{Data Availability} \label{sec:Data}
The developed C++ UMAT subroutine, the complete training dataset (including the densest sampling spectrum), the neural network trained on it, and the Abaqus input file for the perforated plate benchmark are all available under an open-source license at https://doi.org/10.5281/zenodo.18522059. This repository provides the essential components to conduct finite element simulations with the surrogate model and facilitates the adoption and extension of the proposed data-driven framework.

\section{Funding}
We greatefully aknowledge the funding provided by the Deutsche Forschungsgemeinschaft (DFG, German Research Foundation) - project number 414265976 - TRR 277 in the project  project "Bridging Scales – From Geometric Part Details to Construction Elements (C01)" which is part of the collaborative research centre “Additive Manufacturing in Construction - The Challenge of Large Scale”.

\begin{appendices}
\section*{Appendix}
\appendix

\section{FFT-based homogenization}\label{appendix:fft}
\begin{figure}
    \centering
    \input{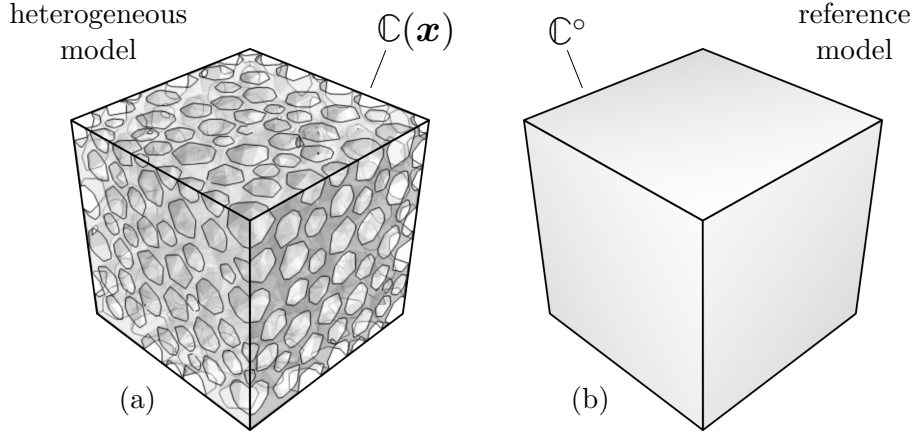}
    \caption{Conceptual basis of mean-field homogenization via the Eshelby formalism: (a) the actual heterogeneous medium and (b) the corresponding auxiliary problem on a homogeneous reference medium.}
    \label{fig:EshelbyFormalism}
\end{figure}

In the spirit of Eshelby's pioneering solution for an ellipsoidal elastic inclusion~\cite{Eshelby1957}, we establish an equivalence between the heterogeneous and homogeneous RVEs shown in \Cref{fig:EshelbyFormalism}. It states that the stress field in the heterogeneous RVE with a spatially varying stiffness tensor $\mathbb{C}(\bm{x})$ can be represented as the stress field of a homogeneous linear elastic reference material, superimposed with a stress polarization field $\bm\tau(\bm{x})$ that arises from the material mismatch. It reads in mathematical terms as 
\begin{equation}
    \bm\sigma(\bm{x}) = \mathbb{C}^\circ:\bm\varepsilon(\bm{x}) + \bm\tau(\bm{x}),
\end{equation}
where $\bm\varepsilon(\bm{x})$ is a spatially varying strain field and $\mathbb{C}^\circ$ is the stiffness tensor of the reference material, which can be chosen arbitrarily without affecting the final outcome, yet this choice has a direct effect on the convergence~\cite{Michel2001}. As a rough guideline, a stiffness tensor representing an average of the actual RVE would be a good choice.

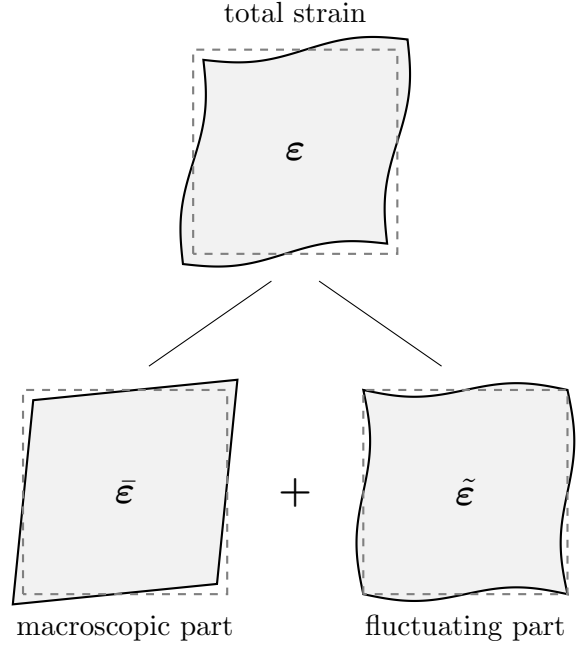
\begin{figure}
    \centering
    \begin{tikzpicture}
    \begin{scope}[scale=0.9]
        \begin{scope}[shift={(4,5)}]
            \filldraw[fill=black!5, draw=black, thick] 
                plot[domain=0:3, samples=100] (\x, {0.1*\x - 0.1*sin(\x*360/3)}) --  
                plot[domain=0:3, samples=100] ({3+0.1*\x-0.1*sin(\x*360/3)}, 0.3+\x) -- 
                plot[domain=3:0, samples=100] (\x+0.3, {3+0.1*\x-0.1*sin(\x*360/3)}) --  
                plot[domain=3:0, samples=100] ({0.1*\x-0.1*sin(\x*360/3)}, \x) -- cycle; 
            \draw[dashed, gray, thick] (0.15,0.15) rectangle (3.15,3.15);
            \node[font=\normalsize, align=center] at (1.65,3.70) {total strain};
            \node[font=\Large, align=center] at (1.65,1.65) {${\bm\varepsilon}$};
        \end{scope}
        \begin{scope}[shift={(1.5,0)}]
            \filldraw[fill=black!5, draw=black, thick] 
                plot[domain=0:3, samples=100] (\x, {0.1*\x - 0.0*sin(\x*360/3)}) --  
                plot[domain=0:3, samples=100] ({3+0.1*\x-0.0*sin(\x*360/3)}, 0.3+\x) -- 
                plot[domain=3:0, samples=100] (\x+0.3, {3+0.1*\x-0.0*sin(\x*360/3)}) --  
                plot[domain=3:0, samples=100] ({0.1*\x-0.0*sin(\x*360/3)}, \x) -- cycle; 
            \draw[dashed, gray, thick] (0.15,0.15) rectangle (3.15,3.15);
            \node[font=\normalsize, align=center] at (1.65,-0.35) {macroscopic part};
            \node[font=\Large, align=center] at (1.65,1.65) {$\bar{\bm\varepsilon}$};
        \end{scope}
        \begin{scope}[shift={(6.5,0)}]
            \filldraw[fill=black!5, draw=black, thick] 
                plot[domain=0:3, samples=100] (0.15+\x, {0.15+0.0*\x - 0.1*sin(\x*360/3)}) --  
                plot[domain=0:3, samples=100] ({3.15+0.0*\x-0.1*sin(\x*360/3)}, 0.15+\x) -- 
                plot[domain=3:0, samples=100] (\x+0.15, {3.15+0.0*\x-0.1*sin(\x*360/3)}) --  
                plot[domain=3:0, samples=100] ({0.15+0.0*\x-0.1*sin(\x*360/3)}, 0.15+\x) -- cycle; 
            \draw[dashed, gray, thick] (0.15,0.15) rectangle (3.15,3.15);
            \node[font=\normalsize, align=center] at (1.65,-0.35) {fluctuating part};
            \node[font=\Large, align=center] at (1.65,1.65) {$\tilde{\bm\varepsilon}$};
        \end{scope}
        \node[font=\Large] at (5.65,1.65) {$\bm{+}$};
        \draw (5.3,4.75) -- (3.5,3.5);
        \draw (6.0,4.75) -- (7.8,3.5);
    \end{scope}
\end{tikzpicture}
    \caption{Additive strain decomposition: the total strain field $\bm\varepsilon$ is decomposed into the superposition of a spatially constant macroscopic part $\bar{\bm\varepsilon}$ and a zero-mean fluctuating field $\tilde{\bm\varepsilon}$.}
    \label{fig:AdditiveSplit}
\end{figure}

The stress tensor $\bm\sigma(\bm{x})$ solves the balance of linear momentum for the actual heterogeneous RVE, which, upon substitution of the above formulation, becomes
\begin{equation}
    \nabla\cdot\bm\sigma(\bm{x}) = \nabla\cdot\bigl(\mathbb{C}^\circ:\bm\varepsilon(\bm{x}) + \bm\tau(\bm{x})\bigr) = \bm{0}.
\end{equation}
Note that body forces are neglected here, as they are accounted for at the macro scale~\cite{Gierden2022}. Next, we adopt the additive strain decomposition illustrated in \Cref{fig:AdditiveSplit}, so that
\begin{equation}
    \bm\varepsilon(\bm{x}) = \bar{\bm\varepsilon} + \tilde{\bm\varepsilon}(\bm{x}),
\end{equation}
where $\bar{\bm\varepsilon}$ is a spatially constant macroscopic field and $\tilde{\bm\varepsilon}(\bm{x})$ is a zero-mean fluctuating field. Since the macroscopic part $\bar{\bm\varepsilon}$ is constant, its divergence vanishes and the balance of linear momentum leads to the governing equation for the fluctuation field
\begin{equation}
    \nabla\cdot\mathbb{C}^\circ:\tilde{\bm\varepsilon}(\bm{x}) = -\nabla\cdot\bm\tau(\bm{x}),
\end{equation}
with the divergence of the stress polarization acting as a source term. Note that solving this equation requires an iterative solution even for linear elastic problems, as the source term itself depends on the unknown fluctuating strain field. This non-homogeneous partial differential equation, complemented with periodic boundary conditions, gives rise to the implicit integral equation
\begin{equation}
    \bm\varepsilon(\bm{x}) = \bar{\bm\varepsilon} - (\mathbb\Gamma^\circ \ast \bm\tau)(\bm{x}),
\end{equation}
where $\mathbb\Gamma^\circ$ is known as the Lippmann--Schwinger operator. This operator has an explicit form in Fourier space that depends only on the properties of the homogeneous linear elastic reference material through
\begin{equation}
    \hat{\mathbb\Gamma}^\circ_{ijkl}(\bm\xi) = \frac{1}{4\mu^\circ|\bm\xi|^2}(\delta_{ik}\xi_j\xi_l + \delta_{il}\xi_j\xi_k + \delta_{jk}\xi_i\xi_l + \delta_{jl}\xi_i\xi_k) - \frac{\lambda^\circ + \mu^\circ}{\lambda^\circ + 2\mu^\circ}\frac{\xi_i\xi_j\xi_k\xi_l}{\mu^\circ|\bm\xi|^4},
\end{equation}
where $\bm\xi$ is the spatial frequency in the Fourier space and $\lambda^\circ$ and $\mu^\circ$ are the Lam\'{e} coefficients of the reference material~\cite{Lucarini2021}. The existence of this explicit form provides the fundamental motivation for introducing a homogeneous reference medium, an idea rooted in Eshelby’s work. Recasting the above equation in Fourier space is a natural step as it is fully consistent with the periodic boundary conditions of the RVE. In addition, the reduction of the convolution integral to a simple pointwise multiplication in the Fourier domain is a beneficial consequence of this transformation.

Note that the convergence of the iterative solver in standard FFT-based homogenization deteriorates with increasing material contrast. It becomes particularly problematic for porous media and breaks down entirely in the limit of infinite contrast. The root of this instability lies in the trigonometric polynomials of the Fourier transform, where the global, edge-to-edge support of the basis functions causes local errors at pore-solid interfaces to propagate throughout the entire domain~\cite{Schneider2016}. A remedy is to rely on finite-difference discretizations, as they have been shown to preserve the well-posedness of the problem and avoid the ill-conditioning arising from high material contrast~\cite{Schneider2021}. One of the most pronounced works belongs to Willot~\cite{Willot2015}, whose methodology has since become a standard discretization scheme in the field~\cite{Lucarini2021}. The framework reformulates the Lippmann--Schwinger operator using centered differences on a rotated staggered grid. The present work adopts this well-established approach. For a detailed derivation and the expression of the modified Lippmann--Schwinger operator, readers are referred to the original work~\cite{Willot2015}.

\Cref{alg:Procedure} summarizes the required steps for solving the periodic boundary value problem on the RVE within the FFT-based homogenization framework. Convergence is reached when the normalized residual of the equilibrium equation falls below a prescribed tolerance, namely
\begin{equation}
\frac{\lVert \nabla \cdot \bm\sigma(\bm{x}) \rVert}{\lVert \bar{\bm\sigma} \rVert} < \epsilon_{\mathrm{tol}},
\end{equation}
where $\bar{\bm\sigma}$ is the homogenized stress passed to the macro-scale upon convergence. It is computed as the volume average of the stress tensor over the RVE domain $\Omega$ as follows:
\begin{equation}
    \bar{\bm\sigma} = \frac{1}{\text{vol}(\Omega)}\int_\Omega \bm\sigma(\bm{x}) \, d\bm{x}.
\end{equation}
This convergence criterion can be efficiently evaluated in Fourier space as
\begin{equation}
\frac{\sum_{\bm\xi \neq \bm{0}} \left\lVert \bm{\xi} \cdot \hat{\bm\sigma}(\bm{\xi}) \right\rVert^2}
{\left\lVert \hat{\bm\sigma}(\bm{0}) \right\rVert^2}
< \epsilon_{\mathrm{tol}}^2,
\end{equation}
where $\hat{\bm\sigma}(\bm{\xi})$ denotes the Fourier transform of the stress field~\cite{Lucarini2021}. This expression follows from two basic properties of the Fourier transform. The numerator arises from the fact that, in Fourier space, the divergence operator becomes a dot product with the frequency vector. This transformation also involves multiplication by the imaginary unit $i$, which, however, does not affect the norm and can therefore be omitted. The denominator, on the other hand, follows from the fact that the zero-frequency component of a field corresponds to its constant (non-fluctuating) part, i.e., its volume average.

\begin{algorithm}
    \begin{algorithmic}[1]
    \State \textbf{Input:} macroscopic strain $\bar{\bm{\varepsilon}}$
    \State \textbf{Output:} homogenized stress $\bar{\bm{\sigma}}$
    \State \textbf{Initialize:} $\bm{\varepsilon}_0(\bm{x}) \gets \bar{\bm{\varepsilon}}$, $k \gets 0$
    \Repeat
        \State $\bm{\tau}_k(\bm{x}) \gets \bm{\sigma}_k(\bm{x}) - \mathbb{C}^\circ:\bm{\varepsilon}_k(\bm{x})$
        \State $\hat{\bm{\tau}}_k(\bm{\xi}) \gets \mathcal{F}\{\bm{\tau}_k(\bm{x})\}$
        \State $\hat{\bm{\varepsilon}}_{k+1}(\bm{\xi}) \gets -\mathbb\Gamma^\circ(\bm{\xi}) : \hat{\bm{\tau}}_k(\bm{\xi})$
        \State $\tilde{\bm{\varepsilon}}_{k+1}(\bm{x}) \gets \mathcal{F}^{-1}\{\hat{\bm{\varepsilon}}_{k+1}(\bm{\xi})\}$
        \State $\bm{\varepsilon}_{k+1}(\bm{x}) \gets \bar{\bm{\varepsilon}} + \tilde{\bm{\varepsilon}}_{k+1}(\bm{x})$
        \State $k \gets k + 1$
    \Until{convergence}
\end{algorithmic}

    \caption{FFT-based solution of the RVE periodic boundary value problem}
    \label{alg:Procedure}
\end{algorithm}

It is worth mentioning that the solution procedure is fundamentally a fixed-point iteration. To significantly accelerate convergence, we implement the limited-memory BFGS algorithm~\cite{Liu1989}. This quasi-Newton method improves the asymptotic convergence rate from linear to superlinear, which dramatically reduces the number of iterations required to reach the solution. Despite its efficiency, multiscale computations, especially for nonlinear materials, are still too computationally involved because \Cref{alg:Procedure} has to be executed at each Gauss point of the global model in every iteration. 

\end{appendices}
\bibliographystyle{ieeetr}
\bibliography{refs,refs_leon}
\end{document}
\endinput